\documentclass[12pt]{amsart}
\usepackage{a4wide, amssymb, bbm, calligra, mathrsfs, stmaryrd}
\usepackage[2cell,arrow,curve,matrix]{xy}
\UseTwocells

\def\calli#1{\textup{\!\textcalligra{#1}\,}}

\let\cal\mathscr

\let\aa\alpha
\let\ee\epsilon
\let\x\times
\let\ho\simeq
\let\Ph\varPhi
\let\ph\varphi
\let\d\partial

\let\bu\bullet
\let\bb\beta
\let\ge\geqslant
\let\le\leqslant
\let\eps\varepsilon
\let\ox\otimes

\let\os\oplus
\let\ot\leftarrow
\let\then\Rightarrow
\let\incl\hookrightarrow

\let\toto\rightrightarrows

\let\al\alpha

\def\op{^{\mathrm{op}}}
\def\mod{\textrm{-}\mathbf{mod}}
\def\1{^{-1}}

\def\PP{{\mathbf P}}
\def\ff{{\boldsymbol f}}
\def\gg{{\boldsymbol g}}
\def\hh{{\boldsymbol h}}
\def\xto#1{\xrightarrow[]{#1}}
\def\xot#1{\xleftarrow[]{#1}}
\DeclareMathOperator\id{id} \DeclareMathOperator\Ob{Ob}
\DeclareMathOperator\End{End} \DeclareMathOperator\Aut{Aut}
\DeclareMathOperator\Mor{Mor} \DeclareMathOperator\Add{Add}
\DeclareMathOperator\Ch{Ch} \DeclareMathOperator\ch{ch}
 \DeclareMathOperator\Der{Der}
\DeclareMathOperator\Ext{Ext} \DeclareMathOperator\Funct{Funct}
\DeclareMathOperator\gldim{gl. dim} \DeclareMathOperator\Hom{Hom}
 \DeclareMathOperator\Inert{Inert}
 \DeclareMathOperator\Ker{Ker}
\DeclareMathOperator\Linext{Linext} \DeclareMathOperator\Nat{Nat}
\DeclareMathOperator\Str{Str} 
\DeclareMathOperator\Tracks{Trext} \DeclareMathOperator\repr{h}

\def\alignbox#1{\begin{aligned}#1\end{aligned}}
\def\set#1{\left\{#1\right\}}
\def\hog#1{\left\llbracket#1\right\rrbracket}

\def\brk#1{\left\langle#1\right\rangle}
\def\alignbox#1{\begin{aligned}#1\end{aligned}}

\newtheorem{Pro}{Proposition}[subsection]
\newtheorem{Le}[Pro]{Lemma}
\newtheorem{The}[Pro]{Theorem}
\newtheorem{Co}[Pro]{Corollary}
\theoremstyle{definition}
\newtheorem{De}[Pro]{Definition}
\theoremstyle{remark}
\newtheorem{Exm}[Pro]{Example}

\def\ext{{\calli{Ext}}}
\def\hom{{\calli{Hom}}}
\def\ta{{\cal T}}

\def\fact{{\mathbf F}}

\def\A{{\mathbb A}}
\def\B{{\mathbb B}}

\def\E{{\mathbb E}}
\def\F{{\mathbb F}}

\def\M{{\mathbb M}}

\def\S{{\mathbb S}}

\def\Z{{\mathbb Z}}

\def\I{{\boldsymbol I}}
\def\bF{{\boldsymbol F}}
\def\bC{{\boldsymbol C}}
\def\bK{{\boldsymbol K}}
\def\bE{{\boldsymbol E}}

\def\a{{\mathscr A}}
\def\b{{\mathscr B}}
\def\c{{\mathscr C}}

\def\f{{\mathscr F}}

\def\o{{\mathscr O}}
\def\r{{\mathscr R}}
\def\s{{\mathscr S}}
\def\u{{\mathscr U}}
\def\g{{\boldsymbol G}}
\def\V{{\mathbf V}}
\def\m{{\boldsymbol m}}
\def\n{{\boldsymbol n}}

\def\CAT{{\mathbf{CAT}}}
\def\Ab{{\mathbbm{Ab}}}
\def\Gr{{\mathbbm{Gr}}}
\def\Fam{{\mathbbm{Fam}}}
\def\AB{{\calli{Ab}\ }}
\def\Set{{\calli{Ens}}}
\def\Theories{{\calli{Theories}}}
\def\Ringoids{{\calli{Ringoids}}}
\def\Th{{\calli{Th}}}
\def\Rings{{\calli{Rings}}}

\def\Monoids{{\calli{Monoids}}}

\begin{document}

\title{Strengthening track theories}
\author{H.-J. Baues}
\address{
Max-Planck-Institut f\"ur Mathematik\\
Vivatsgasse 7\\
Bonn 53111\\
Germany}
\email{baues@mpim-bonn.mpg.de}
\author{M. Jibladze}
\address{
A. Razmadze Mathematical Institute\\
M. Alexidze st. 1\\
Tbilisi 0193\\
Georgia}
\email{jib@rmi.acnet.ge}
\author{T. Pirashvili}
\address{
A. Razmadze Mathematical Institute\\
M. Alexidze st. 1\\
Tbilisi 0193\\
Georgia}
\email{pira@rmi.acnet.ge}

\maketitle

\section{Introduction}

Recent work of the first author on the homotopy category of
4-dimensional manifolds \cite{B4} and on the secondary cohomology
operations \cite{Bsec} is based on the ``calculus of tracks''. One
of the main tricks in \cite{Bsec} is to make some track theories
strong. The aim of this and some subsequent papers is to shed more
light on this procedure. In this paper we prove that certain track
theories are equivalent to strong ones. To be more precise, let us
fix some terminology.

A \emph{track category} is a category enriched in groupoids. Thus
it consists of objects, 1-arrows between them, and 2-arrows, or
tracks, between 1-arrows with the same source and target, and for
two objects $X$, $Y$ of a track category $\ta$ there is their
\emph{Hom-groupoid} $\hog{X,Y}_\ta$, or just $\hog{X,Y}$, whose
objects are 1-arrows $X\to Y$ and morphisms are 2-arrows. Objects
$X$, $Y$ of a track category will be called \emph{homotopy
equivalent} if there are 1-arrows $f:X\to Y$, $g:Y\to X$ with the
composites $fg$, $gf$ isomorphic to identities. A track category
is \emph{abelian} if for any 1-arrow $f:X\to Y$, the group
$\Aut(f)$ of tracks from $f$ to itself is abelian.

Two track categories $\ta$, $\ta'$
are called \emph{weakly equivalent} if there is an enriched functor $F:\ta\to\ta'$
which induces equivalences of hom-groupoids $\hog{X,Y}_\ta\to\hog{FX,FY}_{\ta'}$
and is \emph{essentially surjective}, i.~e. any object of $\ta'$ is
homotopy equivalent to one of the form $FX$.

A \emph{track theory} for us is a track category $\ta$ possessing
finite lax products; this means that for any objects $X$, $Y$ of
$\ta$ there is an object $X\times Y$ with 1-arrows $X\times Y\to
X$, $X\times Y\to Y$ such that the induced functors between
groupoids
$$
\hog{Z,X\times Y}\to\hog{Z,X}\times\hog{Z,Y}
$$
are equivalences of groupoids for all objects $Z$.

A track theory is \emph{strong} if the above functors are in fact
isomorphisms of groupoids.

Morphisms of track theories are enriched functors which are
compatible with lax products. An equivalence of track theories is
a track theory morphism which is a weak equivalence and two track
theories are called \emph{equivalent} if they are made so by the
smallest equivalence relation generated by these.

Our main theorem is that any abelian track theory $\ta$  is
equivalent to a strong one. The fact itself is a trivial
consequence of our results on cohomological properties of
algebraic theories. We believe there exists another, more direct
proof of this, and probably more general result. However the
cohomological results that we obtain are of independent interest
in view of applications to topological Hochschild cohomology
\cite{JP2}, \cite{PW}, \cite{Schwede}.

In  \cite{JP2} the second and third author defined the cohomology
of algebraic theories with some coefficients. In the present paper
we extend this definition in two directions. First, we pass from
single sorted theories to multisorted ones, to obtain our main
theorem in full generality. Second, we extend coefficients for
cohomology. This is necessary for proving our main theorem even
for the particular case of single sorted theories.

For a theory $\A$ we introduce an abelian category $\f(\A)$ in
such a way that the Ext groups in this category yield cohomology
groups of $\A$. This category is in general bigger than the one
introduced in \cite{JP2}, although it is the same in the important
particular case when $\A$ is the theory of modules over a ring ---
see \cite{JPN}. The new cohomologies, just as the old ones, are
closely related to the Baues-Wirsching cohomologies \cite{BW} of
categories; moreover whereas the old coefficients correspond to
the Baues-Wirsching cohomologies of categories with coefficients
in bifunctors, our new extended coefficients correspond to the
Baues-Wirsching cohomologies with coefficients in more general
natural systems.

To get a hint of what new coefficients are, and what kind of
cohomology groups can arise, let us take an example when $\A$ is
the theory of groups $\Gr$. The corresponding coefficient systems
according to \cite{JP2} were functors from the category of
finitely generated free groups to the category of abelian groups.
As we said, in the present paper we consider more general
coefficients, they form the category $\f(\Gr)$, which consists of
assignments $M$ of an $F$-module $M_F$ to each finitely generated
free group $F$, in a way which is functorial in $F$. Then
coefficients in the sense of \cite{JP2} correspond to those
objects $M$ of $\f(\Gr)$ for which the $F$-module structure on
$M_F$ is trivial for all $F$. One typical object of $\f(\Gr)$ is,
for example, $\Omega^1$, which assigns to the group $F$ the
augmentation ideal $\Omega_F^1\subset\Z[F]$ of its group ring,
considered as an $F$-submodule of $\Z[F]$. We will then have, for
any other object $M$ of $\f(\Gr)$, the groups
$\Ext^*_{\f(\Gr)}(\Omega^1,M)$ which will be cohomology groups of
$\Gr$ with coefficients in $M$. We will see that an object similar
to $\Omega^1$ exists in general, and this construction will be
naturally extended to any theory in place of $\Gr$.

One of our main results is that this new cohomology is trivial in
dimensions $>1$ for free theories, just as the old one with more
restricted coefficients. This result together with relationship
between third cohomology group and track extensions \cite{P1},
\cite{P2}, \cite{BD} gives our main result on strengthening of
track theories.


\section{Abelian track categories and cohomology of small categories}\label{prelim}


\subsection{Groupoids, tracks and track
categories}\label{trackdefs} Recall that a \emph{groupoid} is a
category all of whose morphisms are invertible. We will use the
following notation. For a groupoid $\g$, the set of its objects
will be denoted by $\Ob(\g)$ and the set of morphisms by
$\Mor(\g)$. We have the canonical source and target maps
$$
\xymatrix { \Mor(\g)\ar@<1ex>[r]^s\ar@<-1ex>[r]_t&\Ob(\g) }.
$$
A groupoid is called \emph{abelian} if the automorphism group of each
object is an abelian group.

A \emph{2-category} is a category enriched in the category of
small categories. In other words a 2-category $\ta$ consists of  a
class of objects $\Ob(\ta)$, a collection of small categories
$\hog{A,B}=\hog{A,B}_\ta$ for $A,B\in\Ob\ta$ called
hom-\emph{categories} of $\ta$, identities $1_A\in\Ob(\hog{A,A})$
and composition functors $\hog{B,C}\x\hog{A,B}\to\hog{A,C}$
satisfying the usual equations of associativity and identity
morphisms. Objects of the hom-category $f\in\Ob(\hog{A,B})$ are
called \emph{1-arrows} of $\ta$, while morphisms from $\hog{A,B}$
are called \emph{2-arrows}. We will use the following notation for
2-categories. If $f:A\to B$ and $x:B\to C$ are 1-arrows, then the
composite of $f$ and $x$ is denoted by $xf:A\to C$. Notation
$\alpha:f\then f_1$ will indicate a 2-arrow from $f$ to $f_1$,
with $f,f_1\in\Ob(\hog{A,B})$, $A,B\in\Ob(\ta)$. For the
composition of 2-arrows we use additive notation: the identity
2-arrow $f\then f$ of a 1-arrow $f$ will be denoted by $0_f$, and
for 1-arrows $f,g,h:A\to B$ and 2-arrows $\aa:f\then g$,
$\bb:g\then h$, the composite of $\aa$ and $\bb$ in the category
$\hog{A,B}$ is denoted by $\bb+\aa$.

There are several categories associated with a 2-category $\ta$.
The category $\ta_0$ has the  same objects as $\ta$, while
morphisms in $\ta_0$ are 1-arrows of $\ta$. The category $\ta_1$
has the same objects as $\ta_0$. The morphisms $A\to B$ in $\ta_1$
are 2-arrows $\alpha:f\then f_1$ where $f,f_1:A\to B$ are 1-arrows
in $\ta$. The composition in $\ta_1$ is given by $(\bb:x\then
x_1)(\aa:f\then f_1):=(\bb\aa:xf\then x_1f_1)$, where
$$
\bb\aa=\bb f_1+x\aa=x_1\aa+\bb f.
$$
One furthermore has the source and target functors
$$
\xymatrix{
\ta_1\ar@<1ex>[r]^s\ar@<-1ex>[r]_t&\ta_0
},
$$
where $s(\alpha:f\then f_1)=f$ and $t(\alpha:f\then f_1)=f_1$, the
``identity'' functor $i:\ta_0\to\ta_1$ assigning to an 1-arrow $f$
the triple $0_f:f\then f$. Moreover, consider the pullback diagram
$$
\xymatrix{
\ta_1\x_{\ta_0}\ta_1\ar[r]^-{p_2}\ar[d]_{p_1}&\ta_1\ar[d]^t\\
\ta_1\ar[r]^s&\ta_0 };
$$
there is also the ``composition'' functor
$m:\ta_1\x_{\ta_0}\ta_1\to\ta_1$ sending $(\aa:f\then
f_1,\aa':f_2\then f)$ to $\aa+\aa':f_2\then f_1$. Note that these
functors satisfy the identities $sp_1=tp_2$, $sm=sp_2$, $tm=tp_1$
and $si=ti=\id_{\ta_0}$. Sometimes we will also simply write
$\ta_1\toto\ta_0$ to indicate a 2-category $\ta$.

A \emph{track category} $\ta$ is a \emph{category enriched in
groupoids}, i.~e.  is the same as a 2-category all of whose
2-arrows are invertible. If the groupoids $\hog{A,B}$ are abelian
for all $A,B\in\Ob\ta$, then $\ta$ is called an \emph{abelian}
track category. For track categories we might occasionally talk
about \emph{maps} instead of 1-arrows and \emph{homotopies} or
\emph{tracks} instead of 2-arrows. If there is a homotopy
$\alpha:f\then g$ between maps $f,g\in\Ob(\hog{A,B})$, we will say
that $f$ and $g$ are homotopic and write $f\ho g$. Since the
homotopy relation is a natural equivalence relation on morphisms
of $\ta_0$, it determines the \emph{homotopy category}
$\ta_\ho=\ta_0/\ho$. Objects of $\ta_\ho$ are once again objects
in $\Ob(\ta)$, while morphisms of $\ta_\ho$ are homotopy classes
of morphisms in $\ta_0$. For objects $A$ and $B$ we let $[A,B]$
denote the set of morphisms from $A$ to $B$ in the category
$\ta_\ho$. Thus
$$
[A,B]=\hog{A,B}/\ho.
$$
Usually we let $q:\ta_0\to\ta_\ho$ denote the quotient functor. Sometimes
for a 1-arrow $f$ in $\ta$ we will denote $q(f)$ by $[f]$.
A map $f:A\to B$ is a \emph{homotopy equivalence} if
there exists a map $g:B\to A$ and tracks $fg\ho1$ and $gf\ho1$. This is the
case if and only if $q(f)$ is an isomorphism in the
homotopy category $\ta_\ho$. In this case $A$ and $B$ are called \emph{homotopy
equivalent} objects.

A \emph{track functor} $F:\ta\to\ta'$ between track categories is
a groupoid enriched functor. Thus $F$ assigns to each
$A\in\Ob(\ta)$ an object $F(A)\in\Ob(\ta')$, to each map $f:A\to
B$ in $\ta$ --- a map $F(f):F(A)\to F(B)$, and to each track
$\alpha:f\then g$ for $f,g:A\to B$, a track $F(\alpha):F(f)\then
F(g)$ in a functorial way, i.~e. so that one gets functors
$$
F_{A,B}:\hog{A,B}_\ta\to\hog{F(A),F(B)}_{\ta'}.
$$
Moreover these assignments are compatible with identities and
composition, or equivalently induce a functor $\ta_1\to\ta'_1$,
that is, $F(1_A)=1_{F(A)}$ for $A\in\Ob(\ta)$, $F(fg)=F(f)F(g)$,
and $F(\alpha\beta)=F(\alpha)F(\beta)$ for any $\alpha:f\then
f_1$, $\beta:g\then g_1$, $f,f_1:B\to C$, $g,g_1:A\to B$ in $\ta$.

A track functor $F:\ta\to\ta'$ is called a \emph{weak equivalence}
between track categories if the functors
$\hog{A,B}\to\hog{F(A),F(B)}$ are equivalences of groupoids for
all objects $A$, $B$ of $\ta$, and each object $A'$ of $\ta'$ is
homotopy equivalent to some object of the form $F(A)$. Such a weak
equivalence induces a functor $F:\ta_\ho\to\ta'_\ho$ between
homotopy categories which is an equivalence of categories.

\subsection{Preliminaries on cohomology of small categories}\label{factorizdefs}

For us is a crucial fact that any abelian track category defines
an element in the third cohomology group of the corresponding
homotopy category with coefficients in a natural system \cite{P1},
\cite{P2}, \cite{BD}. Therefore we recall the corresponding
notions.

Let $\bC$ be a category. Then the category $\fact\bC$ of
factorizations in $\bC$ is defined as follows. Objects of
$\fact\bC$ are morphisms $f:A\to B$ in $\bC$ and morphisms
$(a,b):f\to g$ in $\fact\bC$ are commutative diagrams
$$
\xymatrix{
A\ar[d]^f&A'\ar[d]^g\ar[l]_a\\
B\ar[r]_b & B' }
$$
in the category $\bC$. A \emph{natural system} on $\bC$  is a
functor $D:\fact\bC\to\AB$ to the category of abelian groups. We
write $D(f)=D_f$. If $a:C\to D$, $f:A\to C$ and $g:D\to B$ are
morphisms in $\bC$, then the induced homomorphism
$(1_A,a)_*:D_f\to D_{af}$ will be denoted by $\xi\mapsto a\xi$,
for $\xi\in D_f$, while $(a,1_B)_*:D_g\to D_{ga}$ will be denoted
by $\eta\mapsto\eta a$, $\eta\in D_g$. We denote by $C^*(\bC;D)$
the following cochain complex:
$$
C^n(\bC;D)=\prod_{\left(A_0\xot{a_1}A_1\xot{\
}\cdots\xot{a_n}A_n\right)\ \in\bC}D_{a_1...a_n},
$$
with the coboundary map given by
\begin{multline*}
d(\ph)(a_1,a_2,...,a_{n+1})=a_1\ph(a_2,...,a_{n+1})+\\
+\sum_{i=1}^n(-1)^if(a_1,...,a_ia_{i+1},...,a_{n+1})+(-1)^{n+1}\ph(a_1,...,a_n)a_{n+1}.
\end{multline*}
According to \cite{BW} the cohomology $H^*(\bC;D)$ of $\bC$ with
coefficients in $D$ is defined as the homology of the cochain
complex $C^*(\bC;D)$.

A morphism of natural systems is just a natural transformation.
For a functor $q:\bC'\to\bC$, any natural system $D$ on $\bC$
gives a natural system $D\circ(\fact q)$ on $\bC'$ which we will
denote $q^*(D)$. There is  a canonical functor
$\fact\bC\to\bC\op\times\bC$ which assigns the pair $(A,B)$ to
$f:A\to B$. This functor allows one to consider any bifunctor
$D:\bC\op\times\bC\to\AB$ as a natural system.  In what follows
bifunctors are considered as natural systems via this
correspondence. Similarly, one has a projection
$\bC\op\times\bC\to\bC$, which yields the functor $\fact\bC\to\bC$
given by $(a:A\to B)\mapsto B$. This allows us to consider any
functor on $\bC$ as a natural system on $\bC$. In particular one
can talk about cohomology of a category $\bC$ with coefficients in
bifunctors and in functors as well. One easily sees that for a
bifunctor $D:\bC\op\times\bC\to\AB$ the group $H^0(\bC;D)$
coincides with the \emph{end} of the bifunctor $D$ (see
\cite{working}), which consists of all families
$(x_C)_{C\in\Ob\bC}$, where $x_C\in D_{1_C}$, for each
$C\in\Ob\bC$, satisfying the condition $a(x_A)=(x_B)a$ for all
$a:A\to B$. In the case of a functor $F:\bC\to\AB$ the group
$H^0(\bC;F)$ is isomorphic to the limit of the functor $F$ and the
groups $H^*(\bC;F)$ are isomorphic to the higher limits (see
\cite{BW}).

\begin{Le}\label{iniciali}
Let $\bC$ be a small category with an initial object $i$. Then for
any functor $F:\bC\to\AB$, one has
\begin{align*}
H^0(\bC;F)&\cong F(i),\\
H^n(\bC;F)&=0 \textrm{ for } n>0.
\end{align*}
\end{Le}

\begin{proof}
In this case, the evaluation of a functor $F$ at $i$ is isomorphic
to the limit of $F$. Thus lim is an exact functor and therefore
higher limits vanish.
\end{proof}

\begin{Exm}
Let $F,G:\bC\to R\mod$ be two functors to the category of left
$R$-modules, for a ring $R$. One can define the bifunctor
$\hom(F,G):\bC\op\times\bC\to\AB$ by
$$
\hom(F,G)(C,D):=\Hom_R(F(C),G(D)).
$$
Then $H^0(\bC;\hom(F,G))\cong\Hom_{\cal A}(F,G)$, where $\cal A$
is the category of all functors from $\bC$ to  $R\mod$. Moreover,
if $F(C)$ is a projective $R$-module for all $C\in\Ob\bC$, then
there is an isomorphism $H^*(\bC;\hom(F,G))\cong\Ext_{\cal
A}^*(F,G)$ \cite{JP2}. We will need a generalization of these
facts, when $R$ is not a constant ring, but a functor from $\bC$
to the category of rings --- see \ref{rtulia} below. In this case
it is necessary to switch to natural systems instead of
bifunctors.
\end{Exm}

Note that for any $\bC$ there is a canonical isomorphism of
categories $\fact\bC\cong\fact(\bC\op)$, which is identity on
objects. Using this we will identify natural systems on $\bC$ and
on $\bC\op$ everywhere in the sequel.

\subsection{ Track extensions and third cohomology of small categories}

As was discovered in \cite{BJ1} if a track category $\ta$ is
abelian, then one has an additional structure. To describe it
we need more notions \cite{BW}, \cite{BD}.

Let $\b$ be a $2$-category. There is a natural system $\End_\b$ of
monoids on $\b_0$ (i.~e. a functor $\fact\b_0\to\Monoids$) which
assigns to an 1-arrow $f:A\to B$ the monoid of all 2-arrows
$f\then f$ in $\b$. Indeed for $g:B\to B'$, $h:A'\to A$ morphisms
in ${\cal B}_0$ we already defined the induced homomorphisms:
\begin{align*}
(\eps\mapsto g\eps):\Hom_\b(f,f)\to \Hom_\b(gf,gf),\\
(\eps\mapsto\eps h):\Hom_\b(f,f)\to \Hom_\b(fh,fh).
\end{align*}

For a track category $\ta$, clearly $\End_\ta=\Aut_\ta$ takes
values in the category of groups. It turns out that the natural
system $\Aut_\ta$ has an additional structure. To describe it let
us introduce the following definition.

\begin{De}\label{connection}
Consider a track category $\ta$. Since taking the category of
factorizations from \ref{factorizdefs} is evidently functorial, applying
it to constituents of $\ta$ gives the diagram
$$
\xymatrix@!C{**[l] \fact(\ta_1\x_{\ta_0}\ta_1)\ar@/^2ex/[r]^{\fact
p_1}\ar[r]|{\fact m}\ar@/_2ex/[r]_{\fact
p_2}&\fact\ta_1\ar@/^2ex/[r]^{\fact s}\ar@/_2ex/[r]_{\fact
t}&\fact\ta_0\ar[l]|{\fact i} },
$$
where the functors $p_1,m,p_2,s,t,i$ are as in \ref{trackdefs}.

A \emph{$\ta$-natural system} with values in a category $\c$ is a
natural system $D:\fact\ta_0\to\c$ on $\ta_0$ together with a
natural transformation $\nabla:D\circ\fact s\to D\circ\fact t$ such that the
diagrams
$$
\xymatrix{
&D\ar@{=}[dl]\ar@{=}[dr]\\
D\circ\fact s\circ\fact i\ar[rr]^{\nabla\fact i}&&D\circ\fact t\circ\fact i
}
$$
and
$$
\xymatrix{
&D\circ\fact s\circ\fact p_1\ar[dl]_{\nabla\fact p_1}\ar@{=}[r]&D\circ\fact t\circ\fact p_2\\
D\circ\fact t\circ\fact p_1\ar@{=}[dr]&&&D\circ\fact s\circ\fact p_2\ar@{=}[dl]\ar[ul]_{\nabla\fact p_2}\\
&D\circ\fact t\circ\fact m&D\circ\fact s\circ\fact m\ar[l]_{\nabla\fact m}
}
$$
commute.

Unfolding this definition in terms of elements one sees easily
that a $\ta$-natural system is the same as a natural system $D$
together with a family of morphisms
$$
\nabla_\xi:D_f\to D_g
$$
in the category $\c$, one for each track $\xi:f\then g$ in $\ta$,
such that the following conditions are satisfied:
\begin{itemize}
\item[i)] $\nabla_{0_f}=\id_{D_f}$ for all 1-arrows $f$ in $\ta$.
\item[ii)] For $\xi:f\then g$, $\eta:g\then h$ one has
$\nabla_{\eta+\xi}=\nabla_\eta\circ\nabla_\xi$.
\item[iii)]
For a diagram
$$
\xymatrix{
\bullet&\bullet\ar[l]_f&\bullet\ltwocell_g^{g_1}{^\xi}&\bullet\ar[l]_h
}
$$
the following diagram
$$
\xymatrix{
D_{fg}\ar[d]_{\nabla_{f\xi}}&D_g\ar[l]_{f\_}\ar[d]_{\nabla_\xi}\ar[r]^{\_h}&D_{gh}\ar[d]_{\nabla_{\xi h}}\\
D_{fg_1}&D_{g_1}\ar[l]_{f\_}\ar[r]^{\_h}&D_{g_1h}
}
$$
commutes.
\item[iv)]
For a diagram
$$
\xymatrix{
\bullet&\bullet\ltwocell_f^{f_1}{^\xi}&\bullet\ar[l]_g&\bullet\ltwocell_h^{h_1}{^\eta}
}
$$
the diagram
$$
\xymatrix{
D_{fg}\ar[dd]_{\nabla_{\xi g}}&&D_{gh}\ar[dd]^{\nabla_{g\eta}}\\
&D_g\ar[ul]_{f\_}\ar[dl]^{f_1\_}\ar[ur]^{\_h}\ar[dr]_{\_h_1}\\
D_{f_1g}&&D_{gh_1}
}
$$
commutes.
\end{itemize}

A morphism $\Ph:(D,\nabla)\to(D',\nabla')$ of $\ta$-natural systems is a natural
transformation $\Ph$ between the functors $D,D':\fact\ta_0\to\c$, such that
the diagram
$$
\xymatrix{
D\circ\fact s\ar[r]^{\Ph\fact s}\ar[d]_\nabla&D'\circ\fact s\ar[d]^{\nabla'}\\
D\circ\fact t\ar[r]^{\Ph\fact t}&D'\circ\fact t}
$$
commutes. We denote by $\ta$-$\Nat$ the category of $\ta$-natural systems.
\end{De}

Let $G:\ta'\to\ta$ be a track functor. For any $\ta$-natural
system $(D,\nabla)$ one defines a $\ta'$-natural system
$G^*(D,\nabla)=(D\circ\fact G,\nabla G)$, where for $\xi':f'\then
g'$ in $\ta'$, $(\nabla G)_{\xi'}:D_{Gf'}\to D_{Gg'}$ is defined
to be $\nabla_{G\xi'}$. In this way one obtains a functor
$$
G^*:\ta\textrm{-}\Nat\to\ta'\textrm{-}\Nat.
$$

\begin{Exm}\label{aut}
For a track category $\ta$, the group-valued natural system $\Aut_\ta$
is equipped with a canonical structure of a $\ta$-natural system given by
$$
\nabla_\xi(a)=\xi+a-\xi.
$$
\end{Exm}

Let $D$ be a natural system on $\ta_\ho$. Then $q^*D$ is a natural
system on $\ta_0$ given by $(q^*D)_f=D_{q(f)}$. Here
$q:\ta_0\to\ta_\ho$ is the canonical projection. Define the
structure of a $\ta$-natural system  on $q^*D$ by
$\nabla=\id:D\circ\fact q\circ\fact s=D\circ\fact q\circ\fact t$.
In this way one obtains the functor
$q^*:\Nat(\ta_\ho)\to\ta$-$\Nat$. Our Theorem \ref{jibladze}
claims that the functor $q^*$ is a full embedding. Actually we
also identify the essential image of the functor $q^*$. We need
the following definition. A $\ta$-natural system $(D,\nabla)$ is
called \emph{inert} if $\nabla_\eps=\id_f$ for all $\eps:f\then
f$. Inert $\ta$-natural systems form a full subcategory of the
category of $\ta$-natural systems, which is denoted by
$\ta$-$\Inert$. It is clear that the image of the functor $q^*$
lies in $\ta$-$\Inert$. It is also clear that $\Aut_\ta$ equipped
with the canonical $\ta$-natural system structure defined in
Example \ref{aut} is inert if and only if $\ta$ is an abelian
track category.

Let us observe that for any track functor $G:\ta'\to\ta$ restriction of the functor
$G^*:\ta$-$\Nat\to\ta'$-$\Nat$ yields the functor $G^*:\ta$-$\Inert\to\ta'$-$\Inert$.

\begin{The}\label{jibladze}
Let $\ta$ be a track category. Then
$q^*:\Nat(\ta_\ho)\to\ta$-$\Inert$
is an equivalence of categories. Furthermore, for any track functor $G:\ta'\to\ta$ the
diagram
$$
\xymatrix{
\Nat(\ta_\ho)\ar[r]^{q^*}\ar[d]^{G_\ho^*}&\ta\textrm{-}\Inert\ar[d]^{G^*}\\
\Nat(\ta'_\ho)\ar[r]^{{q'}^*}&\ta'\textrm{-}\Inert
}
$$
commutes.
\end{The}

\begin{proof}
Let $E$ and $E'$ be natural systems on $\ta_\ho$ and let
$\Ph:q^*E\to q^*E'$ be a morphism of $\ta$-natural systems. We
claim that if $f$ and $g$ are homotopic maps in $\ta_0$ (and
therefore $qf=qg$), then the homomorphisms $\Ph_f:E_{qf}\to
E'_{qf}$ and $\Ph_g:E_{qg}\to E'_{qg}$ are the same. Indeed, we
can choose a track $\xi:f\then g$. Then we have the following
commutative diagram:
$$
\xymatrix{
(q^*E)_f\ar[r]^{\nabla_\xi}\ar[d]_{\Ph_f}&(q^*E)_g\ar[d]^{\Ph_g}\\
(q^*E')_f\ar[r]_{\nabla'_\xi}&(q^*E')_g }
$$
By definition of the $\ta$-natural system structure on $q^*E$ and
$q^*E'$ the morphisms $\nabla_\xi$ and $\nabla'_\xi$ are the
identity morphisms, hence the claim. This shows that the functor
$q^*$ is full and faithful.

It remains to show that for any inert $\ta$-natural system
$(D,\nabla)$ there exists a natural system $E$ on $\ta_\ho$ and an
isomorphism $\Delta:D\to q^*E$ of $\ta$-natural systems. First of
all one observes that if $\xi,\eta:f\then g$ are tracks, then
$\nabla_\xi=\nabla_\eta:D_f\to D_g$. Indeed, thanks to the
property ii) of Definition \ref{connection} we have
$$
\nabla_\xi=\nabla_{\xi-\eta+\eta}=\nabla_{\xi-\eta}\nabla_\eta=\nabla_\eta,
$$
because $\xi-\eta:g\then g$ and $D$ is inert. Therefore for
$qf=qg$ there is a well defined homomorphism $\nabla_{f,g}:D_f\to
D_g$ induced by any track $f\then g$. Then the relation ii) of
Definition \ref{connection} shows that
$\nabla_{g,h}\nabla_{f,g}=\nabla_{f,h}$ for any composable
1-arrows $f,g,h$. By harmless abuse of notation we will just write
$\nabla$ instead of $\nabla_{f,g}$ in what follows.

Since the functor $q:\ta_0\to\ta_\ho$ is identity on objects and
full, we can choose for any arrow $a$ in
$\ta_\ho$ a map $u(a)$ in $\ta_0$ such that $qu(a)=a$. Moreover
for any map $f$ in $\ta_0$ we can choose a track $\delta(f):f\then
u(qf)$.  Now we put
$$
E_a:=D_{u(a)}\textrm{ and
}\Delta_f:=\nabla=\nabla_{f,u(qf)}=\nabla_{\delta(f)}:D_f\to
D_{u(qf)}=E_{qf}.
$$
For a diagram $\xot c\xot a\xot b$ in the category $\ta_\ho$ we
define the homomorphism $c\_:E_a\to E_{ca}$ to be the following
composite:
$$
E_a=D_{u(a)}\xto{u(c)\_}D_{u(c)u(a)}\xto\nabla D_{u(ca)}=E_{ca}.
$$
Similarly we define the homomorphisms $\_b:E_a\to E_{ab}$ to be the following composites:
$$
E_a=D_{u(a)}\xto{\_u(b)}D_{u(a)u(b)}\xto\nabla D_{u(ab)}=E_{ab}.
$$
It follows from the property iii) of Definition \ref{connection}
that for any diagram $\xot{c_1}\xot c\xot a$ in the category $\ta_\ho$
we have the following commutative diagram:
$$
\xymatrix{
D_{u(a)}\ar[d]_{u(c)\_}\ar[dr]^{c\_}\\
D_{u(c)u(a)}\ar[d]_{u(c_1)\_}\ar[r]^\nabla&D_{u(ca)}\ar[d]_{u(c_1)\_}\ar[dr]^{c_1\_}\\
D_{u(c_1)u(c)u(a)}\ar[r]_\nabla&D_{u(c_1)u(ca)}\ar[r]_\nabla&D_{u(c_1ca)}
}
$$
Thus $c_1(c\_)=\nabla(u(c_1)(u(c)\_))$. On the other hand by
definition we have the commutative diagram:
$$
\xymatrix{
D_{u(a)}\ar[d]_{u(c_1c)\_}\ar[dr]^{(c_1c)\_}\\
D_{u(c_1c)u(a)}\ar[r]^\nabla&D_{u(c_1ca)} }
$$
It follows from the property iv) of Definition \ref{connection}
that one has also the following commutative diagram
$$
\xymatrix
{
D_{u(a)}\ar[r]^-{(u(c_1)u(c))\_}\ar[d]_{u(c_1c)\_}&D_{u(c_1)u(c)u(a)}\ar[d]_\nabla\ar[dl]^\nabla\\
D_{u(c_1c)u(a)}\ar[r]_\nabla&D_{u(c_1ca)} }
$$
Therefore
$$
(c_1c)\_=\nabla(u(c_1c)\_)=\nabla(\nabla(u(c_1)(u(c)\_)))=(c_1(c\_)).
$$
Similarly $\_(b_1b)=(\_b_1)b$ and $E$ is a well-defined natural
system on $\ta_\ho$. It remains to show that $\Delta:D\to q^*E$ is
a natural transformation of functors defined on $\fact\ta_0$. To
this end, one observes that for any composable morphisms $g,f$ in
the category $\ta_0$ we have the following commutative diagram
$$
\xymatrix{
D_f\ar[d]^{g\_}\ar[r]^\nabla&D_{u(qf)}\ar[d]^{g\_}\ar[dr]^{u(qg)\_}\\
D_{gf}\ar[r]^\nabla\ar@/_2ex/[rrr]_\nabla&D_{gu(qf)}\ar[r]^\nabla&D_{u(qg)u(qf)}\ar[r]^\nabla&D_{uq(gf)}
}
$$
This means that the following diagram also commutes:
$$
\xymatrix{
D_f\ar[d]^{g\_}\ar[r]^{\Delta_f}&E_{qf}\ar[d]^{(qg)\_}\\
D_{gf}\ar[r]_{\Delta_{gf}}&E_{q(gf)}
}
$$
Similarly the diagram
$$
\xymatrix{
D_g\ar[d]^{\_f}\ar[r]^{\Delta_g}&E_{qg}\ar[d]^{\_(qf)}\\
D_{gf}\ar[r]_{\Delta_{gf}}&E_{q(gf)}
}
$$
also commutes and therefore $\Delta$ is indeed a natural transformation.
\end{proof}

Now let $\ta$ be an abelian track category, so that $\Aut_\ta$ is
a natural system on $\ta_0$ with values in the category of abelian
groups. According to Example \ref{aut} it is equipped with the
canonical structure of a $\ta$-natural system, which is moreover
inert, because $\ta$ is abelian. Thus one can use Theorem
\ref{jibladze} to conclude that there is a natural system $D$
defined on $\ta_\ho$ and an isomorphism of $\ta$-natural systems
$\tau:\Aut_\ta\to q^*D$ defined on $\ta_0$. Roughly speaking a
linear track extension is a choice of such an isomorphism, which
is unique up to a unique isomorphism, in the following sense: if $(D_1,\tau_1)$
is another pair satisfying the same property, then thanks to
Theorem \ref{jibladze} there is a unique isomorphism $\sigma:D\to
D_1$ making the following diagram commute:
$$
\xymatrix@C=.5em {
q^*D\ar[dr]_\tau\ar[rr]^{q^*\sigma}&&q^*D_1\ar[dl]^{\tau_1}\\
&\Aut_\ta }
$$

\begin{De}\label{zveli}
Let $\bC$ be a small category and let $D:\fact\bC\to\AB$ be a
natural system on $\bC$. A linear track extension of $\bC$ by $D$
denoted by
$$0\to D\to \ta_1 \rightrightarrows \ta_0\to \bC \to 0$$
is a pair $(\ta, \tau)$. Here $\ta$ is an abelian track category
equipped with a functor $q:\ta_0\to\bC$ which is full and identity
on objects. In addition for maps $f,g$ in $\ta_0$ we have
$q(f)=q(g)$ iff $f\ho g$. In other words the functor $q$
identifies $\bC$ with $\ta_\ho$. Furthermore
$\tau:q^*D\to\Aut_\ta$ is an isomorphism of $\ta$-natural systems,
where $\Aut_\ta$ is considered as a $\ta$-natural systems as in
Example \ref{aut}.
\end{De}

Hence by virtue of \ref{jibladze} any abelian track category $\ta$
is a  part of the linear track extension
$$0\to D\to \ta_1 \toto \ta_0\to \ta_\ho \to 0,$$
with a natural system $D$, which is defined uniquely up to a
canonical isomorphism.

\

Let $\bC$ be a small category and let $D$ be a natural system on
$\bC$. Objects of the category $\Tracks(\bC;D)$ are linear track
extensions
$$
0\to D\to\ta_1\toto\ta_0\xto q\bC\to0
$$
and morphisms are track functors $F:\ta\to\ta'$ for which $q'F=q$
and $\sigma'=F\sigma$.

\begin{Le}\label{bololema}
Any morphism between track extensions of a small category $\bC$ by
a natural system $D$ is a weak equivalence.
\end{Le}

\begin{proof}
Consider one such morphism represented by $F:\ta\to\ta'$. First of
all, since $F$ is identity on objects, every object of $\ta'$ is
equivalent --- in fact, equal --- to an object of the form $FX$.
Consider now the induced functors between Hom-groupoids
$F_{X,Y}:\hog{X,Y}_\ta\to\hog{FX,FY}_{\ta'}$. These functors are
essentially surjective on objects since $\ta_\ho=\ta'_\ho=\bC$
implies that for any $f':X\to Y$ in $\ta'$ there is an $f:X\to Y$
in $\ta$ with $qf=q'f'$. Then $q'Ff=qf=q'f'$ implies $Ff$ and $f'$
must be homotopic. Next the $F_{X,Y}$ are all full since $q'F=q$
implies that whenever $Ff$ and $Fg$ are homotopic, $f$ and $g$
must be homotopic too, for any $f,g:X\to Y$ in $\ta$. Finally
$\sigma'=F\sigma$ implies that the group homomorphisms
$\Aut_\ta(f)\to\Aut_{\ta'}(Ff)$ are all isomorphisms. This then
clearly implies that all the $F_{X,Y}$ are equivalences of
groupoids.
\end{proof}

\begin{The}\label{h3}{\rm (\cite{P1}, \cite{BD})}
There is a natural bijection
$$
H^3(\bC;D)\approx\pi_0(\Tracks(\bC;D)).
$$
\end{The}

\qed

Here and in what follows $\pi_0$ denotes the set of connected components of a category.

Let
$$
0\to D\to\ta_1\toto\ta_0\xto p\bC\to0
$$
be a linear track extension of $\bC$ by $D$ and let
$f:\bC'\to\bC$ be a functor. Then one can pull back the track extension
to get a linear track extension
$$
0\to D\to\ta'_1\toto\ta'_0\xto{p'}\bC'\to0
$$
of $\bC'$. We define the track category $f^*\ta=\ta'$, as follows.
The objects of $\ta'$ are the same as those of $\bC'$. We will
denote them by $A'$, $B'$ etc. Maps $A'\to B'$ in $\ta'$ are pairs
$(x,\aa)$, where $x:f(A')\to f(B')$ is a morphism in $\ta_0$,
while $\aa:A'\to B'$ is a morphism in $\bC'$ such that
$p(x)=f(\aa)$. If $(x,\aa)$, $(y,\beta)$ are maps $A'\to B'$ in
$\ta'$, then tracks  $(x,\aa)\then(y,\beta)$ in $\ta'$ exist iff
$\aa=\beta$ in $\bC'$. If this condition holds then we put
$$
\Hom_{\hog{A',B'}_{\ta'}}((x,\aa),(y,\beta))=\Hom_{\hog{A,B}_\ta}(x,y),
$$
where $A=f(A')$ and $B=f(B')$. Since the underlying category
$\ta'_0$ is the pullback of $\ta_0\to\bC$ along the functor
$f:\bC'\to\bC$, we will call this construction a pullback
construction. It is clear that one gets a linear track extension
$$
0\to f^*D\to\ta'_1\toto\ta'_0\xto{p'}\bC'\to0,
$$
where $p'$ is identity on objects and on morphisms is given by $p'(x,\aa)=\aa$.
In particular we get the map
$$
f^*:\pi_0(\Tracks(\bC;D))\to\pi_0(\Tracks(\bC';f^*D)); \ \
\ta\mapsto f^*\ta.
$$
One easily checks that in this way we really get a linear track
extension which corresponds to the map $f^*: H^3(\bC;D)\to
H^3(\bC';f^*D)$.

The proof of Theorem \ref{h3} given in \cite{P2} is based on the
following Theorem \ref{fardobiti}, which is going to be crucial in this paper as well.

Let $p:\bK\to\bC$ be a full functor which is identity on objects.
Let $D:\fact\bC\to\AB$ be a natural system on $\bC$. We denote by
$\Tracks(\bC,\bK;D)$ the subcategory of $\Tracks(\bC;D)$ whose
objects are track categories $\ta$ with $\ta_0=\bK$. Morphisms in
$\Tracks(\bC,\bK;D)$ are track functors $\ta\to\ta'$ which are
identity on arrows and hence induce the identity functor
$\bK=\ta_0\to\ta'_0=\bK$. It is clear that $\Tracks(\bC,\bK;D)$ is
a groupoid. In order to relate the set of components
$\pi_0(\Tracks(\bC,\bK;D))$ of $\Tracks(\bC,\bK;D)$ to cohomology
of small categories we need the following variant of the relative
cohomology groups. In the above circumstances $p^*(D)$ is a
natural system on $\bK$, which we will denote still by $D$. Then
$p$ yields a monomorphism of cochain complexes $C^*(\bC;D)\to
C^*(\bK;D)$. We let $C^*(\bC,\bK;D)$ be the cokernel of this
homomorphism. The \emph{$n$-th relative cohomology group}
$H^n(\bC,\bK;D)$ is defined as the $(n-1)$-st homology group of
the cochain complex $C^*(\bC,\bK;D)$. Then one has an exact
sequence
\begin{multline*}\label{grzfar}
0\to H^0(\bC;D)\to H^0(\bK;D)\to H^1(\bC,\bK;D)\to\cdots\\
\to H^n(\bC;D)\to H^n(\bK;D)\to H^{n+1}(\bC,\bK;D)\to\cdots.
\end{multline*}

\begin{The}\label{fardobiti}
Let $p:\bK\to\bC$ be a full functor which is identity on objects
and let $D:\fact\bC\to\AB$ be a natural system. Then there is a
natural bijection
$$
\pi_0(\Tracks(\bC,\bK;D))\approx H^3(\bC,\bK;D).
$$
\end{The}

\begin{proof}
This is exactly Proposition 3.4 of \cite{P2}.
\end{proof}

Let
$$
0\to D\to\ta_1\toto\ta_0\xto p\bC\to0
$$
be a linear track extension of $\bC$ by $D$. According to Theorem
\ref{h3} and Theorem \ref{fardobiti} it defines two elements
$\Ch(\ta)\in H^3(\bC;D)$ and $\ch(\ta)\in H^3(\bC,\ta_0;D)$. It
follows from the proof given in \cite{P2} that
$\partial(\ch(\ta))=\Ch(\ta)$, where $\partial:
H^3(\bC,\ta_0;D)\to H^3(\bC;D)$ is the connecting homomorphism in
the above exact sequence.

\subsection{Lax functors and track extensions}
If two objects $\ta$ and $\ta'$ of the category $\Tracks(\bC;D)$
lie in the same connected component there is no morphism $\ta\to
\ta'$ in $\Tracks(\bC;D)$ in general, but only a diagram of the
form $\ta\ot\ta''\to\ta'$, with an object $\ta''\in
\Tracks(\bC;D)$. The aim of this section is to show that in the
same circumstances there is always a lax functor from $\ta$ to
$\ta'$.

A lax functor $F$ between 2-categories $\ta\to\ta'$ consists of a
map of objects $F:\Ob(\ta)\to\Ob(\ta')$, a collection of functors
$F_{X,Y}:\hog{X,Y}\to\hog{FX,FY}$, for $X,Y\in\Ob(\ta)$, a family
of 2-arrows $o_X:\id_{FX}\then F(\id_X)$ for each $X\in\Ob(\ta)$,
and a natural family of 2-arrows $a_{f,g}:(Ff)(Fg)\then F(fg)$ for
each composable pair of 1-arrows $(f,g)$ in $\ta$. These are
required to satisfy \emph{coherence conditions} --- the following
diagrams
$$
\alignbox{ \xymatrix@!C=2em{
&(F\id_Y)(Ff)\ar@{=>}[dl]_{a_{\id_Y,f}}\\
F(\id_Y
f)&Ff\ar@{=}[l]&(\id_{FY})(Ff)\ar@{=}[l]\ar@{=>}[ul]_{(o_Y)(Ff)}
}}, \alignbox{ \xymatrix@!C=2em{
&(Ff)(F\id_X)\ar@{=>}[dl]_{a_{f,\id_X}}\\
F(f\id_X)&Ff\ar@{=}[l]&(Ff)\id_{FX}\ar@{=}[l]\ar@{=>}[ul]_{(Ff)(o_X)}
}}
$$
and
$$
\alignbox{ \xymatrix@!C {
&(Ff)F(gh)\ar@{=>}[dl]_{a_{f,gh}}\\
F(fgh)&&(Ff)(Fg)(Fh)\ar@{=>}[ul]_{(Ff)a_{g,h}}\ar@{=>}[dl]^{a_{f,g}(Fh)}\\
&F(fg)(Fh)\ar@{=>}[ul]^{a_{fg,h}} } }
$$
must commute for any 1-arrows $f:X\to Y$, $g:W\to X$, $h:V\to W$
in $\ta$.

Let us also explicitate what naturality of $a_{f,g}$ means: it is
equivalent to the commutativity of the diagrams
$$
\alignbox{ \xymatrix{
(Ff)(Fg)\ar@{=>}[d]_{(F\ph)(Fg)}\ar@{=>}[r]^{a_{f,g}}&F(fg)\ar@{=>}[d]^{F(\ph g)}\\
(Ff')(Fg)\ar@{=>}[r]_{a_{f',g}}&F(f'g) } } \textrm{ and }
\alignbox{ \xymatrix{
(Ff)(Fg)\ar@{=>}[d]_{(Ff)(F\psi)}\ar@{=>}[r]^{a_{f,g}}&F(fg)\ar@{=>}[d]^{F(f\psi)}\\
(Ff)(Fg')\ar@{=>}[r]_{a_{f,g'}}&F(fg') } }
$$
for any 1-arrows $f,f':X\to Y$, $g,g':W\to X$ and any 2-arrows
$\ph:f\to f'$, $\psi:g\to g'$ in $\ta$.

A lax functor for which the 2-arrows $o_X$ and $a_{f,g}$ are all
isomorphisms is called a \emph{pseudofunctor}; thus for track
categories these two notions are equivalent. Furthermore a
pseudofunctor is called \emph{strict} if the $o_X$ and $a_{f,g}$
are in fact identities. So a strict pseudofunctor is the same as a
track functor, i.~e. a functor enriched in the category of
categories.

It is immediate from the definitions that a lax functor
$F:\ta\to\ta'$ induces a functor between homotopy categories
$F_\ho:\ta_\ho\to\ta'_\ho$. Moreover it clearly induces monoid
homomorphisms $F_f:\End(f)\to\End(Ff)$ for each 1-arrow $f$ in
$\ta$. A lax functor $F$ between track categories is called a
\emph{lax equivalence} if the functor $F_\ho$ is an equivalence of
categories, and all the group homomorphisms $F_f$ are
isomorphisms. It is easy to see that any lax equivalence is
\emph{locally fully faithful}, i.~e. the functors
$F_{X,Y}:\hog{X,Y}\to\hog{FX,FY}$ are equivalences of groupoids.

\begin{Pro}
Let $D$ be a natural system on a small category $\bC$, and let
$(\ta,\tau)$, $(\ta',\tau')$ be two linear track extensions of
$\bC$ by $D$. Suppose there exists a lax equivalence
$F:\ta\to\ta'$ which is compatible with the track extension
structure in the sense that the triangles
$$
\xymatrix{
&D_\ff\ar[dl]_{\tau_f}\ar[dr]^{\tau'_{Ff}}\\
\Aut_\ta(f)\ar[rr]^{F_f}&&\Aut_{\ta'}(Ff) }
$$
commute for all $\ff:X\to Y$ in $\c$ and all 1-arrows $f$ in $\ta$
with $[f]=\ff$. Then $\Ch(\ta)=\Ch(\ta')\in H^3(\bC;D)$.
\end{Pro}

\begin{proof}
Let us recall how one constructs the characteristic class
$\Ch(\ta)$. For that, one chooses an 1-arrow $s_\ff\in\ff$ in each
homotopy class of 1-arrows in $\ta_\ho=\bC$, and a track
$s_{\ff,\gg}:s_\ff s_\gg\then s_{\ff\gg}$ for each composable pair
of morphisms in $\bC$. Then a cocycle $t$ representing the class
$\Ch(\ta)$ is defined by assigning to a composable triple
$(\ff,\gg,\hh)$ in $\bC$ the element of $D_{\ff\gg\hh}$ given by
the formula
$$
t(\ff,\gg,\hh)=\tau_{s_{\ff\gg\hh}}\1\left(s_{\ff,\gg\hh}+s_\ff
s_{\gg,\hh}-s_{\ff,\gg}s_\hh-s_{\ff\gg,\hh}\right),
$$
where $\tau_{s_{\ff\gg\hh}}:D_{\ff\gg\hh}\to\Aut(s_{\ff\gg\hh})$
is the isomorphism given by the linear track extension structure
of $\ta$. Diagrammatically, $t(\ff,\gg,\hh)$ is the element of
$D_{\ff\gg\hh}$ which corresponds under $\tau_{s_{\ff\gg\hh}}$ to
the automorphism of $s_{\ff\gg\hh}$ given by the counterclockwise
roundtrip in the diagram
$$
\xymatrix@!C{
&s_\ff s_{\gg\hh}\ar@{=>}[dl]_{s_{\ff,\gg\hh}}\\
s_{\ff\gg\hh}&&s_\ff s_\gg s_\hh\ar@{=>}[ul]_{s_\ff s_{\gg,\hh}}\ar@{=>}[dl]^{s_{\ff,\gg}s_\hh}\\
&s_{\ff\gg}s_\hh\ar@{=>}[ul]^{s_{\ff\gg,\hh}}
}.
$$

Given now a lax equivalence $(F,o,a)$ from $\ta$ to $\ta'$ and a
choice of $s_\ff$, $s_{\ff,\gg}$ for $\ta$ as above, we can make
similar choices for $\ta'$ by defining $s'_\ff=Fs_\ff$ and
determining $s'_{\ff,\gg}$ by the commutative diagrams
$$
\xymatrix{
&F(s_\ff s_\gg)\ar@{=>}[dl]_{F(s_{\ff,\gg})}\\
Fs_{\ff\gg}&&(Fs_\ff)(Fs_\gg)\ar@{=>}[ll]_{s'_{\ff,\gg}}\ar@{=>}[ul]_{a_{s_\ff,s_\gg}}
}.
$$
Thus the value of a cocycle $t'$ for $\Ch(\ta')$ on a triple
$\ff,\gg,\hh$ is determined by the outer roundtrip in the diagram
$$
\xymatrix@!C=3em
{
&&Fs_\ff Fs_{\gg\hh}\ar@{=>}[dl]_{a_{s_\ff,s_{\gg\hh}}}\\
&F(s_\ff s_{\gg\hh})\ar@{=>}[dl]_{Fs_{\ff,\gg\hh}}&&Fs_\ff F(s_\gg s_\hh)\ar@{=>}[ul]_{Fs_\ff Fs_{\gg,\hh}}\ar@{=>}[dl]^{a_{s_\ff,s_\gg s_\hh}}\\
Fs_{\ff\gg\hh}&&F(s_\ff s_\gg s_\hh)\ar@{=>}[ul]^{F(s_\ff s_{\gg,\hh})}\ar@{=>}[dl]_{F(s_{\ff,\gg}s_\hh)}&&Fs_\ff Fs_\gg Fs_\hh\ar@{=>}[ul]_{Fs_\ff a_{s_\gg,s_\hh}}\ar@{=>}[dl]^{a_{s_\ff,s_\gg}Fs_\hh}\\
&F(s_{\ff\gg}s_\hh)\ar@{=>}[ul]^{Fs_{\ff\gg,\hh}}&&F(s_\ff s_\gg)Fs_\hh\ar@{=>}[ul]_{a_{s_\ff s_\gg,s_\hh}}\ar@{=>}[dl]^{Fs_{\ff,\gg}Fs_\hh}\\
&&Fs_{\ff\gg}Fs_\hh\ar@{=>}[ul]^{a_{s_{\ff\gg},s_\hh}}
}.
$$
In this diagram, upper and lower squares commute by naturality of
$a$, and the right square commutes as an instance of the coherence
condition. It thus follows that $t'(\ff,\gg,\hh)$ is given by the
roundtrip of the left square, i.~e.
$$
t'(\ff,\gg,\hh)={\tau'_{Fs_{\ff\gg\hh}}}\1\left(Fs_{\ff,\gg\hh}+F(s_\ff
s_{\gg,\hh})-F(s_{\ff,\gg}s_\hh)-Fs_{\ff\gg,\hh}\right).
$$
Now recall that each $F_{X,Y}:\hog{X,Y}\to\hog{FX,FY}$ is a
functor, hence we can write
$$
t'(\ff,\gg,\hh)={\tau'_{Fs_{\ff\gg\hh}}}\1\left(F\left(s_{\ff,\gg\hh}+s_\ff
s_{\gg,\hh}-s_{\ff,\gg}s_\hh-s_{\ff\gg,\hh}\right)\right).
$$
It then follows from compatibility of $F$ with the linear track
extension structures that we obtained a cocycle $t'$ that actually
coincides with $t$.
\end{proof}

Our next aim is to prove the converse of the above proposition,
namely, that if two linear track extensions have the same
characteristic class, then there is a lax equivalence between
them.

For this, let us define for a track category
$\ta=\left(\ta_1\toto\ta_0\right)$ its \emph{relaxation}, which is
a track category $\tilde\ta=(\tilde\ta_1\toto\tilde\ta_0)$
equipped with a weak equivalence $E^\ta:\tilde\ta\to\ta$ (NB: by
definition, \emph{weak} equivalences are \emph{strict} functors,
as opposed to more general \emph{lax} equivalences). We put
$\tilde\ta_0=\PP(\ta_0)$, the path category of $\ta_0$. Recall
that for a graph $G$ its path category $\PP G$ is the free
category on $G$. Thus objects of $\PP G$ are nodes of $G$, and
morphisms of $\PP G$ are finite composable sequences
$$
\xymatrix{ X_0&X_1\ar[l]_{x_1}&\cdots\ar[l]&X_n\ar[l]_{x_n} },
$$
$n\ge0$, of arrows of $G$, identities being empty sequences and
composition given by concatenation. Thus for a small category
$\bC$ considered as a graph there is a canonical functor
$E^\bC:\PP\bC\to\bC$ which is identity on objects, given by
sending a sequence to its composite in $\bC$ (and a sequence with
$n=0$ to the identity of the corresponding object). In particular
$\tilde\ta_0$ comes equipped with such a canonical functor
$E^{\ta_0}:\tilde\ta_0\to\ta_0$. We then define the track
structure of $\tilde\ta$ by pulling it back from $\ta$ along
$E^{\ta_0}$; that is, we define for two sequences of the form
$$
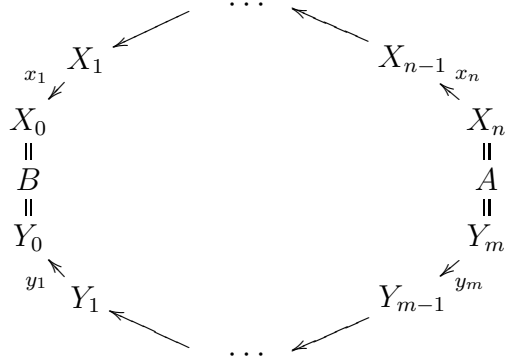
\xymatrix@C=.01em@R=.5em {
                     &                 &&&&&\ \ar[dlllll]&\cdots&&&\\
                     &X_1\ar[dl]_{x_1} &&&&&             &      &&&&&&X_{n-1}\ar[ulllll]\\
X_0                  &                 &&&&&             &      &&&&&&&X_n\ar[ul]_{x_n}\\
B\ar@{=}[u]\ar@{=}[d]&                 &&&&&             &      &&&&&&&A\ar@{=}[u]\ar@{=}[d]\\
Y_0                  &                 &&&&&             &      &&&&&&&Y_m\ar[dl]^{y_m}\\
                     &Y_1\ar[ul]^{y_1} &&&&&             &      &&&&&&Y_{m-1}\ar[dlllll]\\
                     &                 &&&&&\ \ar[ulllll]&\cdots&&&
}
$$
the set of tracks between them by the formula
$$
\Hom_{\hog{A,B}_{\tilde\ta}}((x_1,...,x_n),(y_1,...,y_m))=\Hom_{\hog{A,B}_\ta}(x_1...x_n,y_1...y_m).
$$
Then trivially one checks that this determines a track category,
that $E^{\ta_0}$ extends to a strict functor
$E^\ta:\tilde\ta\to\ta$, that it induces an equivalence (in fact,
an isomorphism) $\tilde\ta_\ho\to\ta_\ho$, and the induced
homomorphisms
$\Aut_{\tilde\ta}((f_1,...,f_n))\to\Aut_\ta(f_1...f_n)$ are all
isomorphisms. In other words, $E^\ta$ is a weak equivalence.

Observe also that there is moreover a canonical lax equivalence
$(F^\ta,o^\ta,a^\ta):\ta\to\tilde\ta$, given as follows:
\begin{itemize}
\item
for an 1-arrow $f:X\to Y$ in $\ta$, $F^\ta_{X,Y}(f)$ is the
1-tuple $(f):X\to Y$ in $\tilde\ta$;
\item for a 2-arrow $\ph:f\then g$, with $f,g:X\to Y$ in $\ta$,
$F^\ta_{X,Y}(\ph)$ is the same 2-arrow
$\ph\in\Hom_{\hog{X,Y}_\ta}(f,g)=\Hom_{\hog{X,Y}_{\tilde\ta}}((f),(g))$;
\item
for $X\in\Ob\bC$, $o^\ta_X:()\then(\id_X)$ is given by
$\id_{\id_X}\in\Hom_{\hog{X,X}_\ta}(\id_X,\id_X)$ $=$
$\Hom_{\hog{X,Y}_{\tilde\ta}}((),(\id_X))$;
\item
for 1-arrows $f:X\to Y$, $g:W\to X$ in $\bC$,
$a^\ta_{f,g}:(f,g)\then(fg)$ is
$\id_{fg}\in\Hom_{\hog{X,Y}_\ta}(fg,fg)$ $=$
$\Hom_{\hog{X,Y}_{\tilde\ta}}((f,g),(fg))$.
\end{itemize}
It is straightforward to check that this indeed defines a lax
equivalence. One notes that the composite of a strict functor and
a pseudofunctor is well defined and in fact $E^\ta F^\ta$ is
identity.

We now have

\begin{The}\label{laxfun}
Let $\ta$, $\ta'$ be linear track extensions of a small category
$\bC$ by a natural system $D$. If $\Ch(\ta)=\Ch(\ta')\in
H^3(\bC;D)$, then there exists a lax equivalence $\ta\to\ta'$.
\end{The}

\begin{proof}
Let us begin by assigning to an 1-arrow $f:X\to Y$ in $\ta$ an
1-arrow $S(f):X\to Y$ in $\ta'$ in such a way that $[S(f)]=[f]:X\to Y$ in
$\bC$. Since $\tilde\ta_0=\PP\ta_0$ is a free category, this
assignment extends uniquely to a functor $S:\tilde\ta_0\to\ta'_0$.
Let us denote by $S^*\ta'$ the track category with
$(S^*\ta')_0=\tilde\ta_0$ obtained by pulling back the 2-arrows
from $\ta'$ along $S$, just as we did when defining $\tilde\ta$.
More precisely, define
$$
\Hom_{\hog{X,Y}_{S^*\ta'}}((x_1,...,x_n),(y_1,...,y_m))
=\Hom_{\hog{X,Y}_{\ta'}}(S(x_1)...S(x_n),S(y_1)...S(y_m)).
$$
Then, exactly as before, one sees that $S$ extends to a strict
functor $S:S^*\ta'\to\ta'$ which is a weak equivalence.

We have now elements $\ch(\ta)\in H^3(\bC,\ta_0;D)$, $\ch(\ta')$
$\in$ $H^3(\bC,\ta'_0;D)$ and $\ch(\tilde\ta)$, $\ch(S^*\ta')$
$\in$ $H^3(\bC,\PP\ta_0;D)$, such that in the diagram of
cohomology groups
$$
\xymatrix{
H^3(\bC,\ta_0;D)\ar[dr]_\d\ar[r]^-{E^*}&H^3(\bC,\PP\ta_0;D)\ar[d]_{\d^{\cong}}&H^3(\bC,\ta'_0;D)\ar[dl]^{\d'}\ar[l]_-{S^*}\\
&H^3(\bC;D) },
$$
one has
\begin{multline*}
\d^{\cong}\ch(\tilde\ta)=\d^{\cong}E^*\ch(\ta)=\d\ch(\ta)=\Ch(\ta)\\
=\Ch(\ta')=\d'\ch(\ta')=\d^{\cong}S^*\ch(\ta')=\d^{\cong}\ch(S^*\ta').
\end{multline*}
Since cohomology of a free category vanishes in dimensions $\ge2$
\cite{BW}, it follows from the long exact sequence connecting the
relative and absolute cohomology groups, that $\d^{\cong}$ is an
isomorphism. Thus $\ch(\tilde\ta)=\ch(S^*\ta')\in
H^3(\bC,\PP\ta_0;D)$. Then there exists an isomorphism of relative
extensions $s:\tilde\ta\to S^*\ta'$, and precomposing it with $S$
we obtain a weak equivalence $Ss:\tilde\ta\to\ta'$. It then
remains to compose this with the lax equivalence
$F^\ta:\ta\to\tilde\ta$ to obtain the required lax equivalence
$SsF^\ta:\ta\to\ta'$.
\end{proof}

\subsection{Linear extensions  and  second
cohomology of categories} To have a more complete picture of the
r\^ole of cohomology of small categories, let us recall the
definition of \emph{linear extensions of categories} and their
relationship with the second cohomology following \cite{BW}. Let
$D$ be a natural system on a small category $\bC$. A linear
extension
$$
0\to D\to\bE\xto p\bC\to0
$$
of $\bC$ by $D$ is a category $\bE$, a full functor $p$ which is
identity on objects, and, moreover, for each morphism $f:A\to B$
in $\bC$, a transitive and effective action of the abelian group
$D_f$ on the subset $p\1(f)\subseteq\Hom_\bE(A,B)$,
$$
D_f\times p\1(f)\to p\1(f); \ \ (a,\tilde f)\mapsto a+\tilde f,
$$
such that the following identity holds
$$
(a+\tilde f)(b+\tilde g)=fb+ag+\tilde f\tilde g.
$$
Here $f$ and $g$ are two composable arrows in $\bC$, $\tilde f\in
p\1(f)$, $\tilde g\in p\1(g)$ and $a\in D_f$, $b\in D_g$. Two linear
extensions $\bE$ and $\bE'$ are \emph{equivalent} if there is an
isomorphism of categories $\ee:\bE\to\bE'$ with $p'\ee=p$ and with
$\ee(a+\tilde f)=a+\ee(\tilde f)$.

Let $\Linext(\bC;D)$ be the set of equivalence classes of linear
extensions of $\bC$ by $D$.

\begin{The}\label{234bw}{\rm(\cite{BW})}
There is a natural bijection
$$
\Linext(\bC;D)\approx H^2(\bC;D).
$$
\end{The}

\section{Local and global $\Ext$-groups}\label{rtulia}

\subsection{The local-global spectral sequence}\label{fmodules}

Let $\I$ be a small category and let $R:\I\to\Rings$ be a functor
to the category of rings with unit. A \emph{left $R$-module} is a
functor $M:\I\to\AB$ together with left $R_i$-module structures on
abelian groups $M_i$, $i\in\Ob(\I)$, such that for any arrow
$\chi:i\to j$ and any $r\in R_i$, $m\in M_i$ one has
$$
R_\chi(rm)=R_\chi(r)M_\chi(m)
$$
in $M_j$. We denote the category of all left modules over a
ring-valued functor $R:\I\to\Rings$ by $R\mod$.

As an example, we can take any small subcategory $\I$ of the
category of commutative rings and let $\o$ be the inclusion
$\I\incl\Rings$. Thus $\o$ is a ring valued functor. For any ring
$S\in\I$ the absolute K\"ahler differentials $\Omega^*_S$ is a
module over $S$. Since $\Omega^*_S$ functorially depends on $S$
we obtain that $\Omega^*\in\o\mod$. Another example comes from
topology. Let $\I$ be a small subcategory of the category of
topological spaces. Then for any ring $R$, the ordinary (singular)
cohomology of spaces with coefficients in $R$ defines a ring
valued functor $H^*(\_;R)$, and for any $R$-module $M$ the functor
$H^*(\_;M)$ is a module over $H^*(\_;R)$ in the above sense.
Similarly $X\mapsto\Z[\pi_1X]$ is a ring valued functor defined on
any small subcategory of the category of pointed topological
spaces, while $X\mapsto\pi_iX$ is a module over it, for any
$i\ge2$.

It is well known that the category $R\mod$ is an abelian category.
Moreover it has enough projective and injective objects (see also
Section \ref{fumodules}). For any $R$-modules $M$ and $N$ one
defines the natural systems $\hom_R(M,N)$ and $\ext^n_R(M,N)$ on
$\I$ by
$$
\hom_R(M,N)_{i\xto\chi j}=\Hom_{R_i}(M_i,N_j)
$$
and
$$
\ext^n_R(M,N)_{i\xto\chi j}=\Ext^n_{R_i}(M_i,N_j)
$$
respectively, where the actions of $R_i$ on $N_j$ are given via
restriction of scalars along $R_\chi:R_i\to R_j$. We call the
natural systems $\hom_R(M,N)$ and $\ext^n_R(M,N)$ \emph{local Hom
and local Ext groups}. One observes that in the case when $R$ is a
constant functor, these natural systems actually come from bifunctors. The
following theorem, which is the main result of this section, was
proved for the particular case of such constant $R$ in \cite{JP2}.

\begin{The}{\bf(the local-to-global spectral sequence)}\label{extlocglobal}
Let $\I$ be a small category and let
$R:\I\to\Rings$ be a functor to the category
of rings with unit. For any $R$-modules $M$ and $N$ there exists a spectral sequence with
$$
E_2^{pq}=H^p(\I;\ext^q_R(M,N))\Longrightarrow\Ext_{R\mod}^{p+q}(M,N).
$$
The result remains true also with rings replaced by ringoids.
\end{The}

The last statement about ringoid-valued functors is essential to
prove our main theorem on strengthening of track theories. We
refer the reader to Section \ref{ringo} for the definition of
ringoids and related stuff, and to page \pageref{prooflocalglobal}
for the proof. Before we go into more detail let us give some
useful consequences.

\begin{Co}
Let $\I$ be a small category and let $M$, $N$ be $\r$-modules,
where
$$
\r:\I\to\Ringoids
$$
is a functor. Then one has a five-term exact sequence
\begin{multline*}
0\to H^1(\I;\hom_\r(M,N))\to\Ext^1_{\r\mod}(M,N)\\
\to H^0(\I;\ext^1_\r(M,N))\to
H^2(\I;\hom_\r(M,N))\to\Ext^2_{\r\mod}(M,N).
\end{multline*}
Moreover, if $\gldim\r_i\le1$ for each object $i$, then
one has an exact sequence
\begin{multline*}
0\to H^1(\I;\hom_\r(M,N))\to\cdots\to H^n(\I;\hom_\r(M,N))\\
\to\Ext^n_{\r\mod}(M,N)\to H^{n-1}(\I;\ext^1_\r(M,N))\to
H^{n+1}(\I;\hom_\r(M,N))\to\cdots
\end{multline*}
\end{Co}

\begin{Co}\label{313proiso}
Suppose $M_i$ is a projective $\r_i$-module for
each $i\in\Ob(\I)$. Then there is an isomorphism
$$
H^*(\I;\hom_\r(M,N))\cong\Ext^*_{\r\mod}(M,N).
$$
\end{Co}

\subsection{Ringoids and modules over them}\label{ringo}

In this subsection we recall some well known facts about ringoids
and modules over them. A good reference on this subject is
\cite{Mitchel}.

A \emph{ringoid} is a category enriched in abelian groups.
It is thus a small category $\r$ together with the structure of
abelian group on its Hom-sets in such a way that composition is biadditive.
Morphisms of ringoids are enriched functors, i.~e. functors preserving the abelian group structures.
These are also called \emph{additive} functors.
The category of ringoids will be denoted by $\Ringoids$.

Let $\r$  be a ringoid. We denote by $\r\mod$ the category of all
covariant additive functors from $\r$ to $\AB$, and by
$\mathbf{mod}$-$\r$ the category of all contravariant additive
functors from $\r$ to $\AB$. Objects from $\r\mod$ are called
\emph{left modules} over $\r$, while those from
$\mathbf{mod}$-$\r$  are called \emph{right modules}.

For any small category $\I$, we let $\Z[\I]$ be the ringoid
with the same objects as $\I$, while for any objects $i$ and
$j$ the group of homomorphisms from $i$ to $j$ in $\Z[\I]$ is
the free abelian group generated by $\Hom_\I(i,j)$:
$$
\Hom_{\Z[\I]}(i,j)=\Z[\Hom_\I(i,j)],
$$
whereas the composition law is induced by
$$
\Z[\Hom_\I(i,j)]\ox\Z[\Hom_\I(j,k)]\cong\Z[\Hom_\I(i,j)\x\Hom_\I(j,k)]\to\Z[\Hom_\I(i,k)].
$$
Then clearly one has $\Z[\I]\mod\simeq\AB^\I$.

For any ringoid $\r$ and an object $c\in\r$ we define
h$_c:\r\to\AB$ and h$^c:\r\op\to\AB$ by
$$
\repr_c(x) = \Hom_\r(c,x)
$$
and
$$
\repr^c(x) = \Hom_\r(x,c).
$$
Then one has natural isomorphisms
$$
\Hom_{\r\mod}(\repr_c,M) \cong M(c)
$$
and
$$
\Hom_{\mathbf{mod}\textrm{-}\r}(\repr^c,N) \cong N(c).
$$

Therefore, the family of objects (h$_c)_{c\in\Ob(\r)}$ (resp. (h$^c)_{c\in\Ob(\r)}$) forms a family of
small projective generators in $\r\mod$  (resp. in $\mathbf{mod}$-$\r$).
The functor h$_c$ is called \emph{the standard free left $\r$-module  concentrated at $c$}.

Let $M:\r\to\AB$ and $N:\r\op\to\AB$ be additive functors. Let
$N\ox_\r M$ be the abelian group defined by
$$
\left.\left(\bigoplus_{c\in\Ob(\r)}N(c)\ox M(c)\right)\right/\sim.
$$
Here $\sim$ is the congruence generated by
$$
N(\alpha)x\ox y \sim x\ox M(\alpha)y
$$
where $\alpha:c_1\to c$ is a morphism in $\r$, $x\in N(c)$, and $y\in M(c_1)$.

Then one has isomorphisms
\begin{align*}
\repr^c\ox_\r M &\cong M(c),\\
N\ox_\r\repr_c &\cong N(c).
\end{align*}

Let $f:\r\to\s$ be a morphism of ringoids.
Composition with $f$ induces a functor
$$
f^*:\s\mod\to\r\mod.
$$
It is well known that $f^*$ has right and left adjoint
functors $f_*$ and $f_!$ respectively (the so-called right and left
Kan extensions) and for any $F:\r\to\AB$ one has isomorphisms
\begin{align*}
(f_*F)(d)&\cong\Hom_{\r\mod}(f^*\repr_d,F),\\
(f_!F)(d)&\cong f^*\repr^d\ox_\r F.
\end{align*}

\subsection{Modules over a ringoid valued functor}\label{fumodules}
Let us consider now a small category $\I$ and a covariant functor
$$
\r:\I\to\Ringoids.
$$
We introduce a category $\smallint_\I\r$ or simply $\smallint\r$ as follows.
Objects of $\smallint\r$ are pairs $(i,x)$, where $i$ is an object of $\I$
and $x$ is an object of  $\r_i$. A morphism $(i,x)\to (j,y)$ is a pair
$(\al,r)$, where $\al:i\to j$ is a morphism in $\I$ and $r:\r_\al(x)\to y$
is a morphism in $\r_j$. Composition in $\smallint\r$ is defined by
$$
(\al,r)\circ(\bb,s)=(\al\circ\bb,r\circ\r_\al(s)).
$$
Then for each $i\in\I$ we have an obvious functor
$\xi_i:\r_i\to\smallint\r$ which assigns $(i,x)$ to an object $x\in\Ob(\r_i)$.

We will say that $M$ is a \emph{left $\r$-module} if the following data are given:
\begin{itemize}
\item[i)] a left $\r_i$-module $M_i$ for each object $i\in\I$;
\item[ii)] a homomorphism $M_\al:M_i\to\r_\al^*M_j$ of $\r_i$-modules for each arrow $\al:i\to j$ of $\I$.
\end{itemize}
Moreover it is required that for any composable morphisms $\al$ and $\bb$ one
has $M_{\al\bb}=M_\al M_\bb$.

If $M$ is a left $\r$-module, $i$ is an object of $\I$, and $x$ is an object of the ringoid
$\r_i$, then we denote by $M_{(i,x)}$ the value $M_i(x)$ of $M_i$ on $x$.
Having this in mind it is clear that a left $\r$-module is nothing else but a functor
$M:\smallint\r\to\AB$ such that each composition $M\circ\xi_i:\r_i\to\AB$, $i\in\I$, is an additive functor.
The category of all left $\r$-modules will be denoted by $\r\mod$.

Yet another description of this category is possible, showing that $\r\mod$ is itself equivalent
to the category of modules over a single ringoid. Given a functor $\r:\I\to\Ringoids$ as above, we define
its \emph{total ringoid} $\r[\I]$ in the following way: the set $\Ob(\r[\I])$ of objects of the ringoid $\r[\I]$
is the disjoint union $\coprod_{i\in\Ob(\I)}\Ob(\r_i)$ --- or else again the set of pairs $(i,x)$,
just as for $\smallint\r$. Morphisms of the ringoid $\r[\I]$ are given by
$$
\Hom_{\r[\I]}((i,x),(j,y))=\bigoplus_{i\xto\alpha j}\Hom_{\r_j}(\r_\alpha(x),y).
$$
Composition homomorphisms are given by
\begin{multline*}
\left(\bigoplus_{i\xto\alpha j}
\Hom_{\r_j}(\r_\alpha(x),y)\right)
\ox
\left(\bigoplus_{j\xto\beta k}\Hom_{\r_k}(\r_\beta(y),z)\right)\\
\xto\cong
\bigoplus_{i\xto\alpha j\xto\beta k}
\Hom_{\r_j}(\r_\alpha(x),y)\ox\Hom_{\r_k}(\r_\beta(y),z)\\
\xto{\bigoplus_{\alpha,\beta}\r_\beta\ox1}
\bigoplus_{i\xto\alpha j\xto\beta k}
\Hom_{\r_k}(\r_\beta\r_\alpha(x),\r_\beta y)\ox\Hom_{\r_k}(\r_\beta(y),z)\\
\xto{\bigoplus_{\alpha,\beta}\circ}
\bigoplus_{i\xto\alpha j\xto\beta k}
\Hom_{\r_k}(\r_\beta\r_\alpha(x),z)
\to
\bigoplus_{i\xto\gamma k}
\Hom_{\r_k}(\r_\gamma(x),z),
\end{multline*}
and the identity of $x\in\Ob(\r_i)$ is the element of
$\bigoplus_{i\xto\eps i}\Hom_{\r_i}(\r_\eps(x),x)$ given by the identity of $x$ in $\r_i$,
situated in the $\id_i$-th summand. It is straightforward to check that this construction
indeed yields a ringoid. One then has

\begin{Pro}\label{total}
For any ringoid-valued functor $\r:\I\to\Ringoids$, the category of left $\r$-modules
is equivalent to $\r[\I]\mod$.
\end{Pro}

\begin{proof}
An $\r[\I]$-module $M$ is a family of abelian groups $(M_{(i,x)})_{x\in\coprod_i\Ob(\r_i)}$
and a family of abelian group homomorphisms
$$
\left(\bigoplus_{i\xto\alpha j}\Hom_{\r_j}(\r_\alpha(x),y)\xto{M_{(i,x),(j,y)}}\Hom_\AB(M_{(i,x)},M_{(j,y)})\right)_{x\in\Ob(\r_i),y\in\Ob(\r_j)},
$$
satisfying certain conditions. Just by universality of sums then, specifying the above homomorphisms
$M_{(i,x),(j,y)}$ is equivalent to specifying families
$$
\left(\Hom_{\r_j}(\r_\alpha(x),y)\xto{M_\alpha}\Hom_\AB(M_{(i,x)},M_{(j,y)})\right)_{\alpha\in\Hom_{\smallint\r}((i,x),(j,y))}.
$$
It is then straightforward to check that the conditions on the $M_{(i,x),(j,y)}$ to form an $\r[\I]$-module
give precisely the conditions on the $M_\alpha$ to form an $\r$-module.
\end{proof}

\begin{Exm}
Since rings with unit are the same as ringoids with single object,
as a particular case of the above construction we have a description of the category
of modules over a ring-valued functor $R:\I\to\Rings$ as in \ref{fmodules}.
Namely, one has an equivalence $R\mod\simeq R[\I]\mod$, where $R[\I]$ is a ringoid
having the same objects as $\I$, with
$$
\Hom_{R[\I]}(i,j)=\bigoplus_{\alpha\in\Hom_\I(i,j)}R_j,
$$
composition given by
$$
\left(\sum_{j\xto\beta k}x_\beta\right)\circ\left(\sum_{i\xto\alpha j}y_\alpha\right)=\sum x_\beta R_\beta(x_\alpha).
$$
To obtain something really familiar, take the further particular case of this,
when $\I$ is a group $G$, considered as a category with one object. Then a ring-valued functor $R$
on this category is the same as a ring $R$ with a $G$-action, and a module $M$ over this functor
is the same as a $G$-equivariant $R$-module, i.~e. an $R$-module with a $G$-action such that
$$
(rm)^g=r^gm^g
$$
for any $r\in R$, $m\in M$, $g\in G$. Furthermore $R[G]$ is in
this case none other than the \emph{crossed group algebra}, i.~e.
the ring obtained by freely adjoining to the multiplicative monoid
of $R$ the group $G$ subject to the commutation relations
$rg=gr^g$ for all $r\in R$, $g\in G$. That $G$-equivariant
$R$-modules are the same as $R[G]$-modules is a classical fact.
\end{Exm}

It is thus clear that $\r\mod$ is an abelian category with enough projective and injective objects.
Let us give the explicit description of the projective generators and injective cogenerators corresponding
to the standard ones from $\r[\I]$.

Take $i\in\Ob(\I)$ and let $x$ be an object of the ringoid $\r_i$.
Then, in accord with the above \ref{total}, associated to the standard free $\r[\I]$-module
concentrated at $(i,x)$ there is a left $\r$-module h$^\r_{i,x}$ given by
$$
\left(\repr^\r_{i,x}\right)_j(y)=\bigoplus_{i\xto\alpha j}\Hom_{\r_j}(\r_\alpha(x),y).
$$
In other words $(\repr^\r_{i,x})_j$ is the direct sum of standard
free $\r_j$-modules:
$$
\left(\repr^\r_{i,x}\right)_j=\bigoplus_{i\xto\alpha j}\repr_{\r_\alpha(x)}.
$$
It follows that for any $\r_j$-module $X$ one has isomorphisms
$$
\Hom_{\r_j}((\repr^\r_{i,x})_j,X)\cong \prod_{i\xto\alpha j}X(\r_\alpha(x)).
$$
Thus for any $\r$-module $M$ one has a natural isomorphism
$$
\Hom_\r(\repr^\r_{i,x},M)\cong M_i(x).
$$

Let now $k$ be an object of $\I$ and let $A$ be an $\r_k$-module.
We denote by $k_*(A)$ the $\r$-module, whose value at $i$ is given
by
$$
(k_*A)_i=\prod_{i\xto\alpha k}\r_\alpha^*A.
$$
The $\alpha$-component of $(k_*A)_i$ has an $\r_i$-module
structure given by restriction of scalars along the ringoid
homomorphism $\r_\alpha:\r_i\to\r_k$. Hence $(k_*A)_i$ is an
$\r_i$-module and now it is clear that $k_*A$ is an $\r$-module.
Moreover the functor $k_*:\r_k\mod\to\r\mod$ is right adjoint to
the evaluation functor ev$_k:\r\mod\to\r_k\mod$, which is given by
ev$_k(M)=M_k$. In particular, if $A$ is an injective $\r_k$-module
then $k_*A$ is an injective $\r$-module. Hence the family
$(k_*Q)_{k,Q}$, is a family of injective cogenerators for the
category of $\r$-modules. Here $k$ runs over the set of objects of
$\I$, and then $Q$ over the set of injective cogenerators of the
category of $\r_k$-modules.

From \ref{total} we also have that for any $M$, $N$ in $\r\mod$ we can calculate their
Ext groups as Ext groups of the corresponding objects in $\r[\I]\mod$. We will denote
$\Ext^*_{\r[\I]}(M,N)$ by $\Ext^*_{\r\mod}(M,N)$ and call them \emph{global Hom and Ext groups}.
One can also define local Hom and Ext functors, exactly as for the ring valued functors.

\begin{Le}\label{Lem:mamuka}
Let us fix $i\in\I$ and $x\in\Ob(\r_i)$.
For any functor $N:\smallint_\I\r\to\AB$ consider the natural system $D$ on $\I$ given by
$$
D_{c\xto\ph d}:=\prod_{i\xto\al c}N(d,\r_{\ph\al}(x)).
$$
Then
$$
H^0(\I;D)=N(i,x).
$$
and
$$
H^n(\I;D)=0 \textrm{ for } n>0.
$$
\end{Le}

\begin{proof}
Consider the \emph{comma category} $i/\I$ (see e.~g. \cite{working});
its objects are arrows $i\to j$, where $j$
runs over objects the category $\I$, and morphism are commutative diagrams
$$
\xymatrix@=1em{
j\ar[rr]&&k\\
\\
&i\ar[luu]\ar[ruu]
}.
$$
One easily checks that
$$
C^*(\I;D)\cong C^*(i/\I;T),
$$
where $T:i/\I\to\AB$ is given by
$$
T\left(i\xto\al c\right)=N(c,\r_\al(x)).
$$
Hence the cohomology of $\I$ with coefficients in $D$ coincides
with the cohomology of the category $i/\I$ with coefficients in
the functor $T$. Since $1_i$ is the initial object in the category
$i/\I$ one can use  Lemma \ref{iniciali} to finish the proof.
\end{proof}

\subsection{Proof of Theorem \ref{extlocglobal}}\label{prooflocalglobal}
We fix a left $\r$-module $N$. We claim that
for any left $\r$-module $X$ one has an isomorphism:
$$
H^0(\I;\hom_\r(X,N))\cong\Hom_{\r\mod}(X,N).
$$
Indeed, it follows from the definition of cohomology that
$H^0(\I;\hom_\r(X,N))$ is isomorphic to the kernel
$$
\Ker\left(
\prod_{i\in\Ob(\I)}\Hom_{\r_i\mod}(X_i,N_i)\to
\prod_{i\xto\alpha j}\Hom_{\r_i\mod}(X_i,N_j)
\right).
$$
Thus $H^0(\I;\hom_\r(X,N))$ consists of families $(f_i:X_i\to
N_i)$ of $R_i$-homomorphisms, such that for any $\alpha:i\to j$
the diagram
$$\xymatrix{
X_i\ar[r]^{f_i}\ar[d]^{X_\alpha}&
N_i\ar[d]^{N_\alpha}\\
X_j\ar[r]^{f_j}&N_j\\
}$$
commutes, and the claim is proved. One observes that the diagram
$$
\xymatrix{
&\Nat(\I)\ar[dr]^{H^0(\I;\_)}\\
\r\mod\op\ar[ur]^{\hom_\r(\_,N)}\ar[rr]_{\Hom_{\r\mod}(\_,N)}&&\AB
}
$$
commutes and the Theorem is a consequence of  the Grothendieck
spectral sequence for composite functors. Of course in order to
apply the Grothendieck theorem  we first have to show that
$H^n(\I;\hom_\r(M,N))=0$ as soon as $n>0$ and $M$ is projective.
To this end we can assume without loss of generality that
$M=h^\r_{i,x}$, for some $i\in\I$ and $x\in\r_i$. In this case
$$
\hom_\r(M,N)_{c\xto\ph d}\cong\prod_{i\xto\al c}N(d,\r_{\ph\al}(x))
$$
and therefore we can use Lemma \ref{Lem:mamuka} to finish the proof.
\qed

\section{Finite product theories}
\subsection{ Basic definitions}\label{theodefs}
A \emph{finite product theory} (simply theory for us) is a small
category with finite products. A morphism of theories is a functor
preserving finite products. With these morphisms, theories form a
category $\Theories$.  Let $\c$ be a category with finite
products. A \emph{model} of a theory $\A$ in the category $\c$,
also termed a \emph{$\c$-valued model of $\A$}, or an
\emph{$\A$-model in $\c$}, is a functor $\A\to\c$ preserving
finite products. Models of $\A$ in $\c$ form a category $\A(\c)$,
with natural transformations as morphisms. Models in the category
$\Set$ of sets will be called simply models, and the category
$\A(\Set)$ will be also denoted by $\A\mod$. It is known that the
category $\A\mod$ is complete and cocomplete for any theory $\A$.
Moreover the inclusion $\A\mod\incl\Funct(\A,\Set)$ preserves all
limits and has a left adjoint, and the Yoneda embedding
$\A\op\to\Funct(\A,\Set)$ factors through it, i.~e. there is a
full embedding $F:\A\op\to\A\mod$. Models in the image of $F$ are
called \emph{finitely generated free models}, so that $\A$ is
equivalent to the opposite of the category of such models. It is
easy to see that the functor $F$ preserves coproducts, i.~e.
$F(X\x Y)$ is a coproduct of $F(X)$ and $F(Y)$ in the category of
models. A morphism of theories $f:\A\to\B$ induces a functor
$$
f^*:\B\mod\to\A\mod,
$$
where $f^*(M)=M\circ f$. Clearly this functor preserves all
limits. Since moreover the categories of models have small
generating subcategories (those of free models), by Freyd's
Special Adjoint Functor Theorem the functor $f^*$ has a left
adjoint
$$
f_!:\A\mod\to\B\mod.
$$
One can see that the square
$$
\xymatrix{
\A\op\ar[r]^-{I_\A}\ar[d]_{f\op}&\A\mod\ar[d]^{f_!}\\
\B\op\ar[r]^-{I_\B}             &\B\mod }
$$
commutes. See \cite{BaWe} for details.

\subsubsection{Single sorted theories}
Let $\S\op\incl\Set$ be the full subcategory of $\Set$ with the
objects $\n=\set{1,...,n}$ for $n\ge0$. Since the category $\S\op$
has finite coproducts, the category $\S$, opposite of the category
$\S\op$ is a theory, which is called the \emph{theory of sets}. To
distinguish objects of $\S$ and $\S\op$ we redenote objects of
$\S$ by $X^0=1$, $X^1=X$, $X^2,X^3,\cdots$. For any $1\le i\le n$
we denote by $x_i:X^n\to X$ the morphism of $\S$ corresponding to
the map $\set{1}\to\n$, which takes $1$ to $i$. It is clear that
$\n$ is a coproduct of $n$ copies of $\set{1}$ in $\S\op$. It follows that
$x_1,...,x_n:X^n\to X$ is a product diagram in $\S$.
One observes that $\S(\c)$ is equivalent to $\c$ for
any category with finite products $\c$. In particular $\S\mod$ is
equivalent to the category $\Set$.

\emph{A single  sorted theory} is a theory morphism $\S\to\A$
which is identity on objects. The full subcategory of
$\S/\Theories$ with single sorted theories as objects will be
denoted by $\Th_1$. Thus objects of single sorted theories are
just natural numbers, which are denoted by  $X^0=1, X^1=X$,
$X^2,X^3,\cdots$. There are projections $x_1,...,x_n$ from $X^n$
to $X$. If $M$ is a model of a single sorted theory $\A$, then
$M(X)$ is called the \emph{underlying set of $M$}. It is then
equipped with operations $u_M:M(X)^n\to M(X)$ for each element $u$ of
$\Hom_\A(X^n,X)$, satisfying identities prescribed by category
structure of $\A$. By this reason, elements of
$\Hom_\A(X^n,X)$ will be called \emph{$n$-ary operations of}
$\A$. Thus for any theory $\A$, the category $\A\mod$ is a variety
of universal algebras. Conversely, for any variety $\V$, the
opposite of the category of the algebras freely generated by the
sets $\n=\set{1,...,n}$, $n\ge0$, is a single sorted theory, whose
category of models is equivalent to $\V$. For example, theory of
groups can be described as follows. Let $\Gr\op$ be the category
with objects $\n$,  $n\ge0$. A morphism from $\n$ to $\m$ is the
same as a homomorphism from the free group on $\n$ to the free
group on $\m$. Clearly $\Gr\op$ is equivalent to the category of
finitely generated free groups. Hence it has finite coproducts.
Therefore the category $\Gr$, the opposite of $\Gr\op$, is a
theory called the \emph{theory of groups}. There is a unique
morphism of theories $\S\to\Gr$ which is identity on objects. Thus
$\Gr$ is a single sorted theory. One observes that $\Gr(\c)$ is
equivalent to the category of group objects in $\c$ and in
particular $\Gr\mod$ is equivalent to the category of groups.

Similarly there is a full embedding
$$
\calli{Rings}\to\Th_1
$$
assigning to a ring $R$ the theory $\M_R$ of left modules over $R$, which is
defined as follows. Let $\M_R$ be the opposite of the full
subcategory of the category $R\mod$  of left $R$-modules with
objects the finitely generated free modules $0$, $R$, $R^2$, ...,
$R^n$, ...; the evident functor $\S\op\to\M_R$ sending $\n$ to
$R^n$ turns $\M_R$ into a single-sorted theory whose category of
models $\M_R\mod$ is equivalent to $R\mod$. Explicitly, the module
corresponding to a model $M$ is $M(R)$, with addition given by
$$
M(R)\x M(R)=M(R^2)\xto{M(+)}{}M(R)
$$
and action of an $r\in R$ given by $M(\_r):M(R)\to M(R)$, where
$\_r:R\to R$ is the homomorphism of left $R$-modules given by
$x\mapsto xr$. For any category $\c$, the category $\M_R(\c)$ is
equivalent to the category of \emph{internal $R$-modules} in $\c$,
i.~e. internal abelian groups $A$ equipped with a unital ring
homomorphism $R\to\End(A)$. In particular, we have the theory of
abelian groups $\Ab=\M_\Z$ such that the category $\Ab(\c)$ is
equivalent to the category of internal abelian groups in $\c$, for
any category $\c$.

\subsubsection{Multisorted theories}
Let $I$ be a set and consider the category $\S\op/I$ of maps
$\n\to I$ for various sets $\n=\set{1,...,n}$.  Morphisms in
$\S\op/I$ from $\n\to I$ to $\m\to I$ are commutative diagrams of
sets
$$
\xymatrix
@=1em{\n\ar[rr]\ar[ddr]&&\m\ar[ddl]\\
\\
&I }
$$
One easily sees that this category has finite coproducts; for
example, coproduct of $f_1:\n_1\to I$ and $f_2:\n_2\to I$ is
$\binom{f_1}{f_2}:\n_1\sqcup\n_2\to I$. in fact, the set of
objects of $\S\op/I$ can be identified with the free monoid
generated by the set $I$ in such a way that a word $i_1...i_n$
represent the coproduct of the objects $i_\nu:{\boldsymbol1}\to
I$, $\nu=1,...,n$. So any $f:\n\to I$ is the coproduct of the
objects $f(1):\boldsymbol1\to I$, ..., $f(n):\boldsymbol1\to I$ in
$\S/I$. We let $\Fam_I$ be the opposite of the category $\S\op/I$.
Then $\Fam_I$ is a theory called the \emph{theory of $I$-indexed
families}. To distinguish objects of $\Fam_I$ and $\S\op/I$ we
denote the object of $\Fam_I$ corresponding to a map $f:\n\to I$
by $X_f$. Hence an object of $\Fam_I$ has the form $X_{i_1}\x...\x
X_{i_n}$ for a unique $n$-tuple $(i_1,...,i_n)\in I^n$. It is
straightforward to check that the functor
\begin{equation*}\label{fam}
\Fam_I(\c)\to\c^I
\tag{*}
\end{equation*}
which assigns to a model $M:\Fam_I\to\c$ the family $M(X_i)_{i\in I}$
is an equivalence.

For a set $I$, an \emph{$I$-sorted theory} is a theory morphism
$\Fam_I\to\A$ which is identity on objects. The full subcategory
of $\Fam_I/\Theories$ with $I$-sorted theories as objects will be
denoted by $\Th_I$.

Although $I$-sorted theories appear to be of very special kind,
one has

\begin{Pro}\label{allsorted}
For any theory $\A$ there is a set $I$ and an $I$-sorted theory
$\Fam_I\to\tilde\A$ such that the category $\tilde\A$ is
equivalent to $\A$.
\end{Pro}

\begin{proof}
Let $I$ be the set $\Ob(\A)$ of objects of $\A$.
We then are forced to take for the set of objects of $\tilde\A$
the free monoid $\sum_{n\ge0}\Ob(\A)^n$ on $I$.
There is an obvious map from this monoid to the set of objects of $\A$,
$\Pi:\Ob(\tilde\A)\to\Ob(\A)$ which assigns to an $n$-tuple $(X_1,...,X_n)$
of objects of $\A$ its product $X_1\x...\x X_n$ in $\A$.
We then simply define
$$
\Hom_{\tilde\A}((X_1,...,X_n),(Y_1,...,Y_m))=\Hom_\A(\Pi(X_1,...,X_n),\Pi(Y_1,...,Y_m)).
$$
This clearly defines the category $\tilde\A$ with the same objects
as $\Fam_{\Ob(\A)}$ and a functor $\tilde\A\to\A$ which is full
and faithful and surjective on objects, i.~e. it is an
equivalence. Moreover by (\ref{fam}) above, models of
$\Fam_{\Ob(\A)}$ in a category with finite products $\c$ are
families $(C_X)_{X\in\Ob(\A)}$ of objects of $\c$, so the
tautological family $(X)_{X\in\Ob(\A)}$ gives a finite product
preserving functor $\Fam_{\Ob(\A)}\to\A$. It is then obvious that
this functor lifts to a functor $\Fam_{\Ob(\A)}\to\tilde\A$ which
is identity on objects.
\end{proof}

A model of an $I$-sorted theory $\Fam_I\to\A$ is just an
$\A$-model. For such a model $\A\to\c$ in a category $\c$ its
\emph{underlying family} is the object of $\c^I$ corresponding to
the composite $\Fam_I\to\A\to\c$. When safe, we will denote images
of morphisms $\omega:X_{i_1}\x...\x X_{i_n}\to X_i$ of $\A$ under
a model $\A\to\c$ by $\omega$ again. Thus intuitively, models $M$
of an $I$-sorted theory $\Fam_I\to\A$ in categories with finite
products $\c$ are $I$-tuples of objects $(C_i)_{i\in I}$,
$C_i=M(X_i)$, equipped with additional structure, namely various
operations of the form
$$
\omega:C_{i_1}\x...\x C_{i_n}\to C_i
$$
corresponding to morphisms $\omega:X_{i_1}\x...\x X_{i_n}\to X_i$ in $\A$.
These operations must further satisfy various identities
expressing the fact that $M$ is a product preserving functor.
In detail, this amounts to the following:
\begin{itemize}
\item
the morphisms corresponding to the projections
$\pi_1:X_{i_1}\x...\x X_{i_n}\to X_{i_1}$
, ...,
$\pi_n:X_{i_1}\x...\x X_{i_n}\to X_{i_n}$
must be product projections themselves;
\item
for morphisms $\omega:X_{i_1}\x...\x X_{i_n}\to X_i$,
$\omega':X_{i'_1}\x...\x X_{i'_m}\to X_i$
and
$\omega_1:X_{i'_1}\x...\x X_{i'_m}\to X_{i_1}$
, ...,
$\omega_n:X_{i'_1}\x...\x X_{i'_m}\to X_{i_n}$
in $\A$ with $\omega(\omega_1,...,\omega_n)=\omega'$, the diagram
$$
\xymatrix@!=2em
{
&C_{i_1}\x...\x C_{i_n}\ar[dl]_-\omega\\
C_i&
&C_{i'_1}\x...\x C_{i'_m}\ar[ul]_-{(\omega_1,...,\omega_n)}
\ar[ll]^-{\omega'}
}
$$
must commute.
\end{itemize}

The ``substrate'' underlying the structure of an $I$-sorted theory
is a family of sets of the form $(S_{(i_1,...,i_n),i})_{(i_1,...,i_n)\in I^n,i\in I}$
for $n=0,1,...$, namely, the sets $\Hom_\A(X_{i_1}\x...\x X_{i_n},X_i)$.
We thus have a forgetful functor
$$
U:\Th_I\to\prod_{n\ge0}\Set^{I^n\x I}.
$$
It is proved in \cite{BV} that this functor admits a left adjoint $F$.
Theories in the image of this left adjoint are \emph{free} theories.
It is more or less obvious that the adjunction counits $FU\A\to\A$
are all full functors, so that in particular one has

\begin{Pro}\label{412enough}  For any theory $\A$ there exists a
morphism $\F\to\A$ from a free theory to $\A$ which is a full
functor.
\end{Pro}

\qed

Moreover, since every componentwise surjective map in $\prod_{n\ge0}\Set^{I^n\x I}$
admits a section, it follows

\begin{Pro}\label{413split} Let $P:\A\to\F$ be a morphism in $\Th_I$
which is a full functor. If $\F$ is a free theory, then $P$ has a
section, i.~e. there is a morphism $S:\F\to\A$ in $\Th_I$ with
$PS=1$.
\end{Pro}
\qed

There is a functor $\Ringoids\to\Theories$. It assigns to a
ringoid $\r$ the theory $\M_\r$ of $\r$-modules. $\M_\r$ is the
additive category freely generated by $\r$, i.~e. it is an
additive category equipped with a homomorphism of ringoids
$I_\r:\r\to\M_\r$ which has the following universal property: for
any additive category $\mathscr A$, precomposition with $I_\r$
induces an equivalence of categories
$$
\Add(\M_\r,{\mathscr A})\cong\Hom_\Ringoids(\r,{\mathscr A}).
$$
There exists an explicit description of $\M_\r$ as the category of
\emph{matrices} over $\r$: $\M_\r$ can be chosen to be an $\Ob(\r)$-sorted
theory, so that its objects are finite families of objects of $\r$, pictured as
$a_1\oplus...\oplus a_n$, for any $a_1,...,a_n\in\r$, $n\ge0$.
Moreover $\Hom_{\M_\r}(a_1\oplus...\oplus a_n,b_1\oplus...\oplus b_m)$ is defined as
$$
\prod_{\substack{i=1,...,m\\j=1,...,n}}\Hom_\r(a_j,b_i),
$$
with composition defined via matrix multiplication, i.~e.
$(f\circ g)_{ik}=\sum_jf_{ij}g_{jk}$ for $f_{ij}:b_j\to c_i$,
$g_{jk}:a_k\to b_j$.

\subsubsection{Tensor product of theories}
In \cite{BV} one finds another useful construction on theories:

\begin{Pro}  For an $I$-sorted theory $\A$ and a
$J$-sorted one, $\B$, there is an $I\x J$-sorted theory $\A\ox B$,
called the \emph{Kronecker product} of $\A$ and $\B$ such that for
any category with products $\c$ one has an equivalence of
categories
$$
\A(\B(\c))\simeq(\A\ox\B)(\c).
$$
\end{Pro}
\qed

\subsubsection{Integrals and cointegrals}
There is a general form of the constructions from \ref{fumodules}. This is
a variation on the \emph{Grothendieck construction, or integral}, which
we briefly recall.

Suppose given a functor $\bF:\I\to\CAT$ from a small category $\I$ to the category of categories,
denoted $(\ph:i\to j)\mapsto(F_\ph:\bF_i\to\bF_j)$.
Then the Grothendieck construction $\smallint_\I\bF$ of $\bF$ is defined as the lax colimit of $\bF$.
Explicitly, it is a category with objects of the form $(i,X)$, with $i\in\Ob(\I)$
and $X\in\Ob(\bF_i)$; morphisms $(i,X)\to(i',X')$ are defined to be pairs $(\ph,f)$,
with $\ph:i\to i'$ and $f:F_\ph(X)\to X'$. Identity morphism for $(i,X)$ is $(\id_i,\id_X)$,
and composition of $(\ph':i'\to i'',f':F_{\ph'}(X')\to X'')$ with $(\ph,f)$ as above
is defined to be the pair $(\ph'\ph,f'F_{\ph'}(f))$. There is a canonical functor
$P_\bF:\smallint_\I\bF\to\I$ given by projection onto the first coordinate,
i.~e. sending $(i,X)$ to $i$ and $(\ph,f)$ to $\ph$.

We will also need the less known \emph{lax limit}, or \emph{cointegral}
of a functor like $\bF$, which we will denote by $\smallint^\I\bF$
(cf. e.~g. \cite[VI 7]{sgaIV}, \cite[5.2.3]{illusie} or
\cite[I,7.12]{fct}). This is equal to the category of sections
of the functor $P_\bF$ above. Thus its objects can be identified
with pairs of families
$$
\left((X_i)_{i\in\Ob(\I)},(F_\ph(X_i)\xto{f_\ph}X_{i'})_{(i\xto\ph i')\in\Mor(\I)}\right)
$$
satisfying $f_{\id_i}=\id_{X_i}$ for any $i\in\Ob(\I)$ and
$f_{\ph'\ph}=f_{\ph'}F_{\ph'}(f_\ph)$ for any $\ph:i\to i'$, $\ph':i'\to i''$.
Whereas morphisms $(X_*,f_*)\to(Y_*,g_*)$ are families $\Ph_i:X_i\to Y_i$ making all the diagrams
$$
\xymatrix@!C=6em{
F_\ph(X_i)\ar[r]^{F_\ph(\Ph_i)}\ar[d]_{f_\ph}&F_\ph(Y_i)\ar[d]^{g_\ph}\\
X_{i'}\ar[r]^{\Ph_{i'}}&Y_{i'}
}
$$
commute.

We will need the following

\begin{Le}\label{prodlax}
A cointegral of theories is a theory. More precisely, suppose
given a functor $\bF:\I\to\CAT$ such that each category $\bF_i$
has finite products and each functor $F_\ph:\bF_i\to\bF_{i'}$
preserves them. Then both categories $\smallint^\I\bF$ and
$(\smallint^\I(\bF\op))\op$ also have finite products which are
computed componentwise, that is, for example, for two objects
$(f_\ph:X_i\to F_\ph(X_{i'}))_{\ph:i\to i'}$ and $(g_\ph:Y_i\to
F_\ph(Y_{i'}))_{\ph:i\to i'}$ of $(\smallint^\I(\bF\op))\op$ their
product is given by the family of the composites
$$
X_i\x Y_i\xto{f_\ph\x g_\ph}F_\ph(X_{i'})\x F_\ph(Y_{i'})\xto\cong F_\ph(X_{i'}\x Y_{i'}).
$$
\end{Le}

Proof is straightforward. \qed

\

Here is an example when such lax limit appears in our context:

\begin{Pro}\label{ringro}
Let $\r:\I\to\Ringoids$ be a
ringoid valued functor on a small category $\I$. Then the category
$\r\mod$ described in \ref{fumodules} is equivalent to
$(\smallint^{\I\op}\bF)\op$ for the functor $\bF:\I\op\to\CAT$ sending
$i\in\Ob(\I)$ to the category $\r_i\mod\op$ of modules over the
ringoid $\r_i$ and $\ph:i\to j$ --- to the ``restriction of
scalars'' functor $\r_\ph^*:\r_j\mod\to\r_i\mod$ induced by the
ringoid homomorphism $\r_\ph:\r_i\to\r_j$.
\end{Pro}

\begin{proof}
This follows straightforwardly from the definition of the category $\r\mod$
in \ref{fumodules}.
\end{proof}

A construction similar to that of \ref{fumodules} can be performed with theories too.

Let $\A:\I\to\Theories$ be a functor from a small category $\I$ to
theories. Define a \emph{model} of $\A$ to be a collection
$(M_i)_{i\in\Ob(\I)}$ of $\A_i$-models, one for each
$i\in\Ob(\I)$, together with a collection of morphisms
$M_\ph:M_i\to\A_\ph^*M_{i'}$, one for each $\ph:i\to i'$ in
$\Mor(\I)$, such that $M_{\id_i}=\id_{M_i}$ and
$M_{\ph'\ph}=\A_\ph^*(M_{\ph'})M_\ph$ for any $\ph:i\to i'$,
$\ph':i'\to i''$. Thus also straightforwardly one has

\begin{Pro}
For any functor $\A:\I\to\Theories$, the category of $\A$-models
is equivalent to $(\smallint^{\I\op}F)\op$, where $F:\I\to\CAT$ assigns
the category $\A_i\mod$ to $i\in\Ob(\I)$ and the functor $\A_\ph^*$
to $\ph:i\to i'$.
\end{Pro}
\qed

\subsubsection{Comma category as models} As an application of
previous discussion we prove that the comma category of a
category of models of a theory is still a category of models for a
theory.

\begin{Pro}\label{groth}
For an $I$-sorted theory $\A$
and any model $M$ in $\A\mod$, the category $\smallint_\A M$ is a
$\left(\coprod_{i\in I}M_i\right)$-sor\-ted theory and moreover
the comma category $\A\mod/M$ is equivalent to the category
of models $\left(\smallint_\A M\right)\mod$.
\end{Pro}

\begin{proof}
Any object $N$ of $\A\mod$ equipped with a morphism $f:N\to M$
can be considered as a collection of sets
$$
\left(N_x=f_A\1(x)\subseteq N(A)\right)_{x\in\coprod_{A\in\Ob(\A)}M(A)}
$$
and maps $N_{x_1}\x...\x N_{x_n}\to N_{\omega(x_1,...,x_n)}$,
for all $(x_1,...,x_n)\in M(X_{i_1})\x...\x M(X_{i_n})$ and
$\omega:X_{i_1}\x...\x X_{i_n}\to X_i$ in $\A$,
fitting into certain commutative diagrams.

Then regarding $M$ as an object of $\Set^\A$, and defining
$N(x)=N_{M(p_1)x}\x...\x N_{M(p_n)x}$, for $x\in M(X_{i_1}\x...\x
X_{i_n})$, we can consider the above data as a functor $\tilde
N:\smallint_\A M\to\Set$, which sends the object $x\in
M(X_{i_1}\x...\x X_{i_n})$ of the latter category to the product
of the objects $\tilde N(X_{i_\nu})$, $\nu=1,...,n$. Now the proof
follows from the subsequent lemma.
\end{proof}

\begin{Le}
A functor $M:\A\to\Set$ preserves finite products
if and only if the category $\smallint_\A M$ has finite
products and the canonical functor $P:\smallint_\A M\to M$,
sending $m\in M(X)$ to $X$, preserves them.
\end{Le}

\begin{proof}
Let us first recall that
functors of the form $P:\smallint_\A M\to\A$ for any functor
$M:\A\to\Set$ are characterized by a property called \emph{discrete opfibration}:
\begin{quote}
for any $x\in\smallint_\A M$ and any $\ph:Px\to a$, there is a unique $\psi:x\to y$
with $P\psi=\ph$.
\end{quote}
Using this property it is easy to prove that a pullback of a product preserving
discrete fibration between categories with products along a product preserving functor
is again a product preserving functor between categories with products.

The ``only if'' part then follows because of the following pullback diagram
in the category of categories
$$
\xymatrix
{
\smallint_\A M\ar[r]\ar[d]^P&\Set_\bu\ar[d]_U\\
\A\ar[r]^M&\Set
}
$$
in which $\Set_\bu$ denotes the category of pointed sets and $U$
the forgetful functor: since the latter is a discrete opfibration
and preserves products, it follows that $\smallint_\A
M$ will have and $P:\smallint_\A M\to\A$ preserve them too.

For the ``if'' part, we again use the discrete fibration property to prove
\begin{itemize}
\item[a)] $M(1)$ has single element: the particular case of the above discrete opfibration condition
with $Px=a=1$ implies that for any $x\in P\1(1)$ one has
$\left(x\xto{\id_x}{}x\right)$ $=$ $\left(x\xto{!_x}{}1\right)$, since $P(\id_x)=P(!_x)=\id_1$.
\item[b)] $M(a_1\x a_2)\xto{(M\pi_1,M\pi_2)}{}Ma_1\x Ma_2$ is bijective: this follows
from another two particular cases of the discrete opfibration condition --- with $x=x_1\x x_2$ for
some $x_i\in P\1(a_i)$ and $\ph=\pi_i$, $i=1,2$;
indeed these cases give that there are unique $\psi_i$ starting out of $x$ with $P(\psi_i)=\pi_i$,
hence $x$ is a unique element of $M(a_1\x a_2)$ satisfying $M\pi_i(x)=x_i$, $i=1,2$.
\end{itemize}
\end{proof}

\begin{Co}\label{commod}
For any theory $\A$ and any functor $M:\I\to\A\mod$, there is an equivalence
$$
\A\mod^\I/M\simeq\A/M\mod,
$$
where $\A/M:\I\to\Theories$ is the functor given by $i\mapsto\smallint_\A M_i$.
\end{Co}

\begin{proof}
An object of $\A\mod^\I/M$ consists of homomorphisms of $\A$-models $p_i:N_i\to M_i$, $i\in\Ob(\I)$,
and $N_\ph:N_i\to N_{i'}$, $\ph:i\to i'$, such that all squares
$$
\xymatrix{
N_i\ar[r]^{N_\ph}\ar[d]_{p_i}&N_{i'}\ar[d]^{p_{i'}}\\
M_i\ar[r]^{M_\ph}&M_{i'}
}
$$
commute and moreover $N_{\id_i}=\id_{N_i}$, $N_{\ph'\ph}=N_{\ph'}N_\ph$
for all $i$ and all composable pairs $\ph'$, $\ph$.

It is then clear that such data can be equivalently figured out
as a collection of objects $(N_i,p_i)\in\A\mod/M_i$, $i\in\Ob(\I)$,
together with a collection of morphisms $(N_i,p_i)\to M_\ph^*(N_{i'},p_{i'})$,
for $\ph:i\to i'$ in $\I$ which satisfy exactly the conditions determining
an object of $(\smallint^{\I\op}(\A/M(\_))\op)\op$, i.~e., by definition,
of $\A/M\mod$. It is straightforward to check that this correspondence also
carries over to morphisms.
\end{proof}

\subsection{Enveloping ringoids}

\begin{Pro}\label{envel}  For any $I$-sorted theory $\A$ there exists a
ringoid $U(\A)$, depending functorially on $\A$, such that $\Ab(\A\mod)$ is equivalent to the category of
$U(\A)$-modules.
\end{Pro}

\begin{proof}
The key observation here is that in the presence of an abelian group structure
any operation like $\omega:X_1\x...\x X_n\to X$ must be an abelian group homomorphism,
hence have the form $\omega(x_1,...,x_n)=\omega_1(x_1)+...+\omega_n(x_n)$ for some
unary operations $\omega_i:X_i\to X$.

Let the set of objects of $U(\A)$ be $I$, and present morphisms of $U(\A)$
by generators and relations as follows. For each $\omega:X_{i_1}\x...\x X_{i_n}\to X_i$ in $\A$
we pick $n$ generators $\brk{\omega,1}:X_{i_1}\to X_i$, ..., $\brk{\omega,n}:X_{i_n}\to X_i$.
And for each such $\omega$ and any $\omega_1:X_{i'_1}\x...\x X_{i'_m}\to X_{i_1}$, ...,
$\omega_n:X_{i'_1}\x...\x X_{i'_m}\to X_{i_n}$ we impose the relations
$$
\brk{\omega(\omega_1,...,\omega_n),\mu}=\sum_{\nu=1}^n\brk{\omega,\nu}\circ\brk{\omega_\nu,\mu}
$$
for $\mu=1,...,m$. So a $U(\A)$-module is a collection of abelian groups $(A_i)_{i\in I}$ and homomorphisms
$\brk{\omega,\nu}:A_{i_\nu}\to A_i$, $\omega\in\Hom_\A(X_{i_1}\x...\x X_{i_n},X_i)$, $\nu=1,...,n$
satisfying the above relations. Then from any such module we obtain an object of $\Ab(\A\mod)$
by defining
$$
\omega(a_1,...,a_n)=\sum_{\nu=1}^n\brk{\omega,\nu}a_\nu
$$
for $\omega$ as above and $(a_1,...,a_n)\in A_{i_1}\x...\x A_{i_n}$. Conversely, if $(A_i)_{i\in I}$
is given the structure of an object from $\Ab(\A\mod)$, then we define
\begin{align*}
\brk{\omega,\nu}a=\omega(0,...,0,&a,0,...,0).\\[-1ex]
\nu\textrm{-th}&^\uparrow\textrm{position}
\end{align*}

It is easy to see that these procedures determine mutually inverse
equivalences between the category of $U(\A)$-modules and
$\Ab(\A\mod)$.
\end{proof}

One then has

\begin{Co}
For a theory $\A$, the following conditions are equivalent:
\begin{itemize}
\item[i)] $\A$ is isomorphic to $\M_\r$ for some ringoid $\r$;
\item[ii)] $\A\mod$ is an additive category;
\item[iii)] the canonical homomorphism of theories $\A\to\Ab\ox\A$ is an isomorphism;
\item[iv)] $\A$ is isomorphic to $\Ab\ox\B$ for some theory $\B$.
\end{itemize}
\end{Co}

\begin{proof}
Implication i)$\then$ii) is clear, ii)$\iff$iii)
follows from the fact that a category $\a$ with finite products is additive
iff the forgetful functor $\Ab(\a)\to\a$ is an equivalence, and iii)$\then$iv) is trivial.
Finally iv)$\then$i) follows from the above proposition.
\end{proof}

\begin{Co}\label{commab}
For any model $M$ of a theory $\A$, there exists a ringoid
$\u(M)$, the \emph{enveloping ringoid} of $M$, depending functorially on $M$,
such that the
category $\Ab(\A\mod/M)$ is equivalent to the category of
$\u(M)$-modules.
\end{Co}

\begin{proof}
Of course this is just a particular case of the previous
proposition in view of \ref{groth}. Let us, however, give explicit
presentation of $\u(M)=U(\smallint_\A M)$ in this case, assuming
for simplicity that $\A$ is an $I$-sorted theory. The set of
objects of $\u(M)$ is then $\coprod_{i\in I}M(X_i)$, and the
morphisms are generated by ones of the form
$\brk{\omega,x_1,...,x_n,\nu}:x_{i_\nu}\to\omega(x_1,...,x_n)$,
for each $\omega\in\Hom_\A(X_{i_1}\x...\x X_{i_n},X_i)$,
$(x_1,...,x_n)\in
 M(X_{i_1})\x...\x M(X_{i_n})$
and $\nu\in\set{1,...,n}$. The defining relations are indexed by
data $\omega\in\Hom_\A(X_{i_1}\x...\x X_{i_n},X_i)$,
$\omega_1\in\Hom_\A(X_{i'_1}\x...\x X_{i'_m},X_{i_1})$, ...,
$\omega_n\in\Hom_\A(X_{i'_1}\x...\x X_{i'_m},X_{i_n})$,
$(x_1,...,x_m)\in M(X_{i'_1})\x...\x M(X_{i'_m})$, and
$\mu\in\set{1,...,m}$ and have the form
$$
\brk{\omega(\omega_1,...,\omega_n),x_1,...,x_m,\mu}
=\sum_{\nu=1}^n
\brk{\omega,\omega_1(x_1,...,x_m),...,\omega_n(x_1,...,x_m),\nu}
\circ\brk{\omega_\nu,x_1,...,x_m,\mu}.
$$
Once again, functoriality is obvious from this presentation.
\end{proof}

Occasionally we will write $\u_\A(M)$ to make explicit dependence on $\A$. This construction
is known under various names in the literature --- see e.~g. \cite{BE} or \cite{R}.
We will also need a generalization of this fact to functors, which requires the following

\begin{Le}\label{modint}
Given a theory $\A$ and a functor $\bF:\I\to\CAT$, there is an equivalence
$$
\A\left(\smallint^\I\bF\right)\simeq\smallint^\I\A(\bF_{\_}),
$$
where $\A(\bF_{\_}):\I\to\CAT$ is given by $i\mapsto\A(\bF_i)$.
\end{Le}

\begin{proof}
It is easy to see that for any category $\A$ whatsoever there is an equivalence
$$
\Funct(\A,\smallint^\I\bF)\simeq\smallint^\I\bF^\A,
$$
where $\bF^\A:\I\to\CAT$ is given by $i\mapsto\Funct(\A,\bF_i)$.
On the other hand we know by \ref{prodlax}
that products in $\smallint^\I\bF$ are computed componentwise;
this implies easily that an object of $\smallint^\I\bF^\A$
like $f_\ph:M_i\to F_\ph M_i'$, with $M_i:\A\to\bF_i$, etc.
corresponds to an $\A$-model in $\smallint^\I\bF$ iff
each $M_i$ is an $\A$-model in $\bF_i$.
\end{proof}

\begin{Pro}
For any theory $\A$ and any functor $M:\I\to\A\mod$, there is an
equivalence of categories
$$
\Ab(\Funct(\I,\A\mod)/M)\simeq\u(M)\mod,
$$
where $\u(M):\I\to\Ringoids$ is the functor defined by
$i\mapsto\u(M_i)$.
\end{Pro}

\begin{proof}
From \ref{commod}, there is an equivalence
$$
\Ab(\Funct(\I,\A\mod)/M)\simeq\Ab(\A/M\mod);
$$
but by definition
$$
\A/M\mod=\left(\smallint^{\I\op}(\A\mod/M_{\_})\op\right)\op,
$$
Hence by \ref{modint}, there is an equivalence
$$
\Ab(\A/M\mod)=\left(\smallint^{\I\op}\bF\right)\op,
$$
where $\bF:\I\to\CAT$ is the functor given by $i\mapsto\Ab(\A\mod/M_i)$.
Now by \ref{envel}, $\bF$ is isomorphic to the functor $\u(M_{\_})\mod:\I\to\CAT$
given by $i\mapsto\u(M_i)\mod$; and we have proved in \ref{ringro} that
there is an equivalence
$$
\left(\smallint^{\I\op}\u(M_{\_})\mod\right)\op\simeq\u(M)\mod.
$$
\end{proof}

Given a theory $\A$, its model $M\in\A\mod$, and an object $p:A\to
M$ of the category $\Ab(\A\mod/M)\simeq\u_\A(M)\mod$, we will
denote by $\Der(M;A)$ the abelian group of all sections of $A\to
M$, i.~e. the set of all morphisms $s:M\to A$ of $\A$-models with
$ps=1_M$. Elements of $\Der(M;A)$ will be called
\emph{derivations} of $M$ in $A$. $\Der(M;A)$ is contravariantly
functorial in $M$, in the following sense. For a morphism $f:M'\to
M$ of models we get the induced homomorphism
$f^*:\Der(M;A)\to\Der(M';f^*A)$, where $f^*A$ denotes the pullback
of $p:A\to M$ along $f$. Equivalently, one might interpret
$\Der(M';f^*A)$ as the abelian group of all $\A$-model morphisms
$M'\to A$ over $M$, i.~e. fitting in the commutative diagram
\begin{equation}\label{sectri}
\alignbox{
\xymatrix{
&A\ar[d]^p\\
M'\ar[ur]\ar[r]^f&M
}
}.
\tag{\ddag}
\end{equation}
Clearly also $\Der(M;A)$ is covariantly functorial in $A$ and so
defines a functor $\Der(M;\_)$ on $\u_\A(M)\mod$. We then have

\begin{Pro}\label{436de}
The functor $\Der(M;\_)$ is representable. That is, there exists
an $\u_\A(M)$-module $\Omega^1_M$ with natural isomorphism
$\Der(M;A)\cong\Hom_{\u_\A(M)}(\Omega^1_M,A)$ for all $A$.
Moreover $\Omega^1$ depends functorially on $M$. When $M$ is a
finitely generated free $\A$-model, then $\Omega^1_M$ is a
projective object of $\u(M)\mod$.
\end{Pro}

\begin{proof}
Following the equivalence from \ref{commab}, for an $\u(M)$-module $A$ the corresponding
object of $\Ab(\A\mod/M)$ is the $\A$-model with $X_i\mapsto\coprod_{x\in M(X_i)}A(x)$,
with the $\A$-model structure assigning to $\omega:X_{i_1}\x...\x X_{i_n}\to X_i$ the
operation
$$
\omega:\coprod_{x_1\in M(X_{i_1})}A(x_1)\x...\x\coprod_{x_n\in M(X_{i_n})}A(x_n)\to
\coprod_{x\in M(X_i)}A(x)
$$
given by
$$
\omega(a_1,...,a_n)=\sum_{\nu=1}^n\brk{\omega,x_1,...,x_n,\nu}a_\nu.
$$
Then
$$
\Der(M;A)\subset\prod_{\substack{i\in I\\x\in M(X_i)}}A(x)
$$
consists of those families $(d(x)\in A(x))_{x\in\coprod_iM(X_i)}$ which respect all these operations.
That is, $\Der(M;A)$ consists of assignments, to each $x\in M(X_i)$, of an element $d(x)\in A(x)$,
in such a way that for any $\omega:X_{i_1}\x...\x X_{i_n}\to X_i$ and any $x_\nu\in M(X_{i_\nu})$,
$\nu=1,...,n$, one has
\begin{equation*}\label{der}
d(\omega(x_1,...,x_n))=\sum_{\nu=1}^n\brk{\omega,x_1,...,x_n,\nu}d(x_\nu).
\tag{*}
\end{equation*}
Because of this expression it is natural to call such assignments \emph{derivations}.

We then present $\Omega^1_M$ by generators and relations as a $\u(M)$-module as follows:
it has generators $d(x)\in\Omega^1_M(x)$ for each $x\in M(X_i)$ and each $i\in I$;
and the defining relations are (\ref{der}) above. It is then clear that $\Omega^1_M$ carries a generic
derivation $d$, so that one has a natural isomorphism
$$
\Hom_{\u(M)}(\Omega^1_M,A)\xto\cong\Der(M;A)
$$
given by $f\mapsto fd$. That $\Omega^1$ is functorial in $M$ is also clear from the construction.

Now suppose $M$ is a finitely generated free model $F(X)$, i.~e.
there is an $X\in\A$ with $M=\Hom_\A(X,\_)$. Then it is
straightforward to check using Yoneda lemma that for an object of
$\Ab(\A\mod/M)$ corresponding to a $\u(M)$-module $A$ we will have
$\Der(F(X);A)\cong A(\id_X)$. It follows that
$\Hom_{\u(F(X))}(\Omega^1_{F(X)},A)$ is an exact functor of $A$,
i.~e. $\Omega^1_{F(X)}$ is projective. In fact of course this
actually means that $\Omega^1_{F(X)}=\repr_{\id_X}$.
\end{proof}

\section{Cartesian natural systems}

\subsection{Definitions, motivation, examples}

Let $\A$ be a theory and let $D$ be a natural system on $\A$. We
will say that the natural system $D$ is \emph{cartesian} (or
\emph{compatible with products} --- cf. \cite{BTONKS})  if
for any product diagram $p_k:X_1\x...\x X_n\to X_k$, $k=1,...,n$
and any morphism $f:X\to X_1\x...\x X_n$ the homomorphism
$$
D_f\to D_{p_1f}\x...\x D_{p_nf}
$$
given by $a\mapsto(p_1a,...,p_na)$ is an isomorphism. Obviously
$D$ is cartesian if and only if it satisfies the above condition
with $n=0$ and $n=2$, i.~e.
\begin{itemize}
\item $D_{!_X}=0$ for the unique morphism $!_X:X\to1$ to the
terminal object;
\item $D_f\to D_{p_1f}\x D_{p_2f}$ is an isomorphism for any
$f:X\to X_1\x X_2$.
\end{itemize}
One observes that if a bifunctor $D:\A\op\x\A\to\AB$ preserves
products in the second variable, then the natural system induced
by $D$ is cartesian. We denote by $\f(\A)$ the category of
cartesian natural systems on $\A$.

\begin{Exm}\label{dercart}
Recall that in \ref{fumodules} we have defined the notion of a module over a
ringoid-valued functor. Let us then, for a theory $\A$, consider the ringoid-valued
functor $\u_\A$ on $\A$ given by $X\mapsto\u_A(F(X))$, where $F(X)=\Hom_\A(X,\_)$
is the free finitely generated $\A$-model corresponding to the object $X$.
For any two objects $A$, $B$ of $\u_\A\mod$, similarly to the natural systems $\hom$ defined in
\ref{fmodules}, there is a natural system $\hom(A,B)$ on $\A$ given by
$$
\hom(A,B)_{X\xto fY}=\Hom_{\u_\A(F(Y))}(A_Y,F(f)^*B_X),
$$
where the ringoid morphism $F(f):\u_\A(F(Y))\to\u_\A(F(X))$ is
induced by $F(f):F(Y)\to F(X)$, i.~e. by $(g\mapsto
gf):\Hom_\A(Y,\_)\to\Hom_\A(X,\_)$. Let us find out when is this
natural system cartesian. For this it will be convenient to
rewrite the above in the following way:
$$
\hom(A,B)_{X\xto fY}=\Hom_{\u_\A(F(X))}(F(f)_!A_Y,B_X).
$$
Indeed as we saw in \ref{theodefs} all the functors $F(f)^*$ have left adjoints.
The above conditions then show that this natural system is cartesian if and only if
\begin{itemize}
\item $\Hom_{\u(F(X))}(F(!_X)_!A_1,B_X)=0$ for all $X$;
\item the canonical morphism
$$
\Hom_{\u(F(X))}(F(f)_!A_{X_1\x X_2},B_X)\to\Hom_{\u(F(X))}(F(p_1f)_!A_{X_1}\oplus F(p_2f)_!A_{X_2},B_X)
$$
is an isomorphism for any $f:X\to X_1\x X_2$.
\end{itemize}
In particular $\hom(A,B)$ is cartesian for \emph{all} $B$ if and only if $A$ satisfies
\begin{itemize}
\item $A_1=0$;
\item $F(p_1)_!A_{X_1}\oplus F(p_2)_!A_{X_2}\to A_{X_1\x X_2}$ is an isomorphism for any $X_1$, $X_2$.
\end{itemize}
It is natural to call such an $A$ a \emph{cartesian $\u_\A$-module}.

We already have a nice example of such: the $\Omega^1$ constructed above. Indeed
any $\u_\A$-module $B$ determines a natural system $\Der(\_;B)$ on $\A$ in the
following way: for a morphism $f:X\to Y$ of $\A$, put
$$
\Der(\_;B)_f=\Der(F(Y);f^*(B_X)).
$$
Here $p_X:B_X\to F(X)$ is the object of $\Ab(\A\mod/F(X))$ corresponding to $B(X)$
under the equivalence $\u_\A(F(X))\mod\simeq\Ab(\A\mod/F(X))$.
That this is indeed a natural system, follows from the functorial properties of $\Der$.
Moreover this natural system is cartesian. Indeed, $\A$-models of the form $F(X)$ are
the representable ones, $F(X)(Y)=\Hom_\A(X,Y)$. Then considering
the diagram (\ref{sectri}) we see that $\Der(F(Y);f^*(B_X))$ can be identified with the set of all elements
$b\in B_X(Y)$ with $p_X(b)=f\in F(X)(Y)=\Hom_\A(X,Y)$. Then given $f_i:X\to X_i$, $i=1,...,n$,
one has
\begin{align*}
\Der(\_;B)&_{(f_1,...,f_n)}=\Der(F(X_1\x...\x X_n);(f_1,...,f_n)^*(B_X))\\
&\approx\set{b\in B_X(X_1\x...\x X_n)\ \mid\ p_X(b)=(f_1,...,f_n)}\\
&\approx\set{(b_1,...,b_n)\in B_X(X_1)\x...\x B_X(X_n)\ \mid\ p_X(b_i)=f_i, i=1,...,n}\\
&\approx\Der(\_;B)_{f_1}\x...\x\Der(\_;B)_{f_n}.
\end{align*}
But it is immediate from \ref{436de} that there is an $\u_\A$-module $\Omega^1$ such that
the natural system $\Der(\_;B)$ is actually isomorphic to $\hom(\Omega^1_{F(\_)},B)$.
Namely, $\Omega^1$ is just given by $X\mapsto\Omega^1_{F(X)}$. It is then a cartesian $\u_\A$-module,
i.~e. one has
\begin{itemize}
\item $\Omega^1_{F(1)}=0$;
\item $F(p_1)_!\Omega^1_{F(X_1)}\oplus F(p_2)_!\Omega^1_{F(X_2)}\to\Omega^1_{F(X_1\x X_2)}$
is an isomorphism for any $X_1$, $X_2$.
\end{itemize}
\end{Exm}

The following fact goes back to \cite{JP2}.
\begin{Le}\label{compatible}
Let
$$
0\to D\to\E\xto P\A\to0
$$
be a linear extension of a theory $\A$ by a natural system
$D$. Then $D$ is cartesian iff $\E$ is a theory and $P$ is a theory morphism.
\end{Le}

\begin{proof}
Take a product diagram $p_i:X_1\x...\x X_n\to X_i$,
$i=1,...,n$, and choose arbitrarily $\tilde p_i$ in $\E$ with $P(\tilde p_i)=p_i$.
This then gives a commutative diagram
$$
\xymatrix@C=8em{
\Hom_\E(X,X_1\x...\x X_n)\ar[d]_P\ar[r]^-{\tilde f\mapsto(\tilde p_1\tilde f,...,\tilde p_n\tilde f)}&\Hom_\E(X,X_1)\x...\x\Hom_\E(X,X_n)\ar[d]^P\\
\Hom_\A(X,X_1\x...\x X_n)\ar[r]^-\approx&\Hom_\A(X,X_1)\x...\x\Hom_\A(X,X_n)
}
$$
which shows that $\E$ has and $P$ preserves finite products iff all the maps
$$
P\1(f)\to P\1(p_1f)\x...\x P\1(p_nf),
$$
given by $\tilde f\mapsto(\tilde p_1\tilde f,...,\tilde p_n\tilde f)$ are bijective.

On the other hand the above maps are equivariant with respect to the group homomorphisms
$$
D_f\to D_{p_1f}\x...\x D_{p_nf}
$$
and the actions given by the linear extension structure. Our proposition then follows
from the following easy lemma.
\end{proof}

\begin{Le}
Suppose given a group homomorphism $f:G_1\to G_2$ and an $f$-equivariant map $x:X_1\to X_2$
between sets $X_i$ with transitive and effective $G_i$-actions. Then $x$ is bijective
iff $f$ is an isomorphism.
\end{Le}

\begin{proof}
See e.~g. \cite[Lemma 3.5]{JP2}
\end{proof}

\begin{The}
There is an equivalence of categories
$$
\Phi:\f(\A)\to\u_\A\mod;
$$
in particular, $\f(\A)$ is an abelian category with enough projectives and injectives.
Moreover the quasi-inverse of this equivalence assigns to an object $A$
of $\u_\A\mod$ the cartesian natural system $\Der(\_;A)$ from \ref{dercart}.
\end{The}

\begin{proof}
As always, we can assume here that $\A$ is an $I$-sorted theory.
Then for a cartesian natural system $D$ on $\A$, to define
$\Phi(D)$ we must first name for each $X\in\Ob(\A)$ a
$\u_\A(F(X))$-module $\Phi(D)_X$. The set of objects of
$\u_\A(F(X))$ is $\coprod_{i\in I}F(X)(X_i)$ (see \ref{commab}),
i.~e. $\coprod_{i\in I}\Hom_\A(X,X_i)$. We then define values of
$\Phi(D)_X$ on these objects by
$$
\Phi(D)_X(X\xto x X_i)=D_x.
$$
Next action of morphisms of $\u_\A(F(X))$ is uniquely determined by requiring, for
$(x_1,...,x_n):X\to X_{i_1}\x...\x X_{i_n}$ and $\omega:X_{i_1}\x...\x X_{i_n}\to X_i$,
commutativity of the diagrams
$$
\xymatrix{
& D_{(x_1,...,x_n)}\ar[dl]_{\omega\_}& D_{x_1}\x...\x D_{x_n} \ar[l]_-\cong
\\
D_{\omega(x_1,...,x_n)} &&& D_{x_\nu},
\ar[lll]_-{\brk{\omega,x_1,...,x_n,\nu}} \ar[ul]_{\iota_\nu}
}
$$
where the isomorphism is the inverse of the canonical map that is required by cartesianness of $D$,
and $\iota_\nu$ is the $\nu$-th embedding into $\os=\x$ of abelian groups.

We also have to define action on $\Phi(D)$ of morphisms $f:X\to Y$
in $\A$, which must be $\u_\A(F(Y))$-module morphisms
$\Phi(D)_Y\to F(f)^*(\Phi(D)_X)$, where the functor
$F(f)^*:\u_\A(F(X))\mod\to\u_\A(F(Y))\mod$ is the restriction of
scalars along the ringoid morphism $\u_\A(F(Y))\to\u_\A(F(X))$
induced by the morphism of $\A$-models $F(f):F(Y)\to F(X)$. Now
$F(f)^*(\Phi(D)_X)$ is easily seen to be given by $(y:Y\to
X_i)\mapsto D_{yf}$, so what we must choose is a suitably
compatible family of abelian group homomorphisms
$$
\Phi(D)_f(Y\xto y X_i):D_y\to D_{yf},
$$
and these we declare to be the action of $\_f$ on $D$. It is then
straightforward that all of the above indeed gives a functor
$\Phi:\f(\A)\to\u_\A\mod$.

Next note that, as we have seen in \ref{436de}, one has
$\Der(F(X);A)\cong A(\id_X)$ for any $\u_\A(F(X))$-module $A$, so
in particular for any $f:X\to Y$ in $\A$ we have by \ref{dercart}
$$
\Der(\_;\Phi(D))_f=\Der(F(Y);F(f)^*(\Phi(D)_X))\cong
F(f)^*(\Phi(D)_X)(\id_Y)=D_{\id_Yf}=D_f.
$$
Conversely, given a $\u_\A$-module $A$, by definition
\begin{multline*}
\Phi(\Der(\_;A))_X(X\xto xX_i)=\Der(\_;A)_x=\Der(F(X_i);F(x)^*(A_X))\\
\cong F(x)^*(A_X)(\id_{X_i})=A_X(x).
\end{multline*}
(Of course one should also check these on morphisms, but this is
straightforward too).
\end{proof}

\subsection{Cohomology of theories}

For a theory $\A$ and an object $A\in\u_A\mod$, we next define the cohomology
$$
H^*(\A;A)
$$
by the equality
$$
H^*(\A;A):=\Ext_{\u_\A\mod}(\Omega^1,A).
$$
Here $\Omega^1$ is turned into an object of $\u_\A\mod$ as in \ref{dercart} above.

The following is the main theorem of this section

\begin{The}\label{jib}
Let $\A$ be a  theory and let $A$ be an $\u_\A$-module. Then
$$
H^*(\A;A)\cong H^*(\A;\Der(\_;A)),
$$
where on the left we have cohomology just defined, while on the right --- the Baues-Wirsching
cohomology of the category $\A$ with coefficients in the natural system $\Der(\_;A)$.
\end{The}

\begin{proof}
By Proposition \ref{436de} one has an isomorphism of natural
systems:
$$
\Der(\_;A)\cong\hom(\Omega^1,A).
$$
Hence the result is a consequence of Corollary \ref{313proiso}. The fact
that the condition of Corollary \ref{313proiso} holds follows from
Proposition \ref{436de}.
\end{proof}

\begin{Co}
If $\F$ is a free $I$-sorted theory and $D$ is a cartesian natural system on $\F$,
then
$$
H^i(\F;D)=0, \ \ i>1.
$$
\end{Co}

\begin{proof}
First consider the case $i=2$; thanks to Theorem \ref{234bw}
it suffices to show that any linear extension of $\F$ by $D$ splits. By
Lemma \ref{compatible} any such extension is an extension in $\Theories$
and we can use Proposition \ref{413split} to conclude that it really splits.
If $i\ge3$ we can use Theorem \ref{jib} to pass to
the theory cohomologies.
The latter are Ext-groups in appropriate abelian categories vanishing
on injective objects and in dimension two we can use the long cohomological sequence
associated to an extension
$$
0\to A\to I\to B\to 0
$$
with injective $I$ to finish the proof.
\end{proof}

This result in the case when $D$ is a bifunctor over a single sorted theory
was proved in \cite{JP2} (see Proposition 4.22 of loc. cit.).

\section{Track theories and the proof of the main theorem}

\subsection{Track theories}

Recall that
a track category $\ta$ is said to have \emph{finite lax products}, if for any
finite collection $X_1$, ..., $X_n$ of its objects there exists a
family of 1-arrows $p_1:X\to X_1$, ..., $p_n:X\to X_n$ such that the
induced functors
$$
\hog{Y,X}\to\hog{Y,X_1}\x...\x\hog{Y,X_n},
$$
$f\mapsto(p_1f,...,p_nf)$, are equivalences of groupoids for any object $Y$ of $\ta$.
A track category with finite lax products is called a \emph{track theory}.
A track theory is called \emph{strong} if the products in it are in fact strong, i.~e.
the above equivalences of groupoids are in fact isomorphisms.

As we saw a linear extension of a theory $\A$ by a natural system
$D$ is again a theory provided $D$ is cartesian. The situation
changes dramatically for track extensions. Let $\ta$ be a linear
track extension of a theory $\A$ by a cartesian natural system
$D$. Then in general $\ta$ is not a strong track theory, but only
a track theory. This is the subject of the following
\begin{Pro}\label{dekartulitrack}
For a linear track extension
$$
D\to\ta_1\toto\ta_0\to\A,
$$
the corresponding track category is a track theory if and only if
$\A$ is a theory and $D$ is a cartesian natural system.
\end{Pro}

\begin{proof}
Assume $\ta$ is a track theory, so that
$$
\hog{Y,X}\to\hog{Y,X_1}\x...\x\hog{Y,X_n}
$$
are equivalences of groupoids. Thus they induce bijections
on connected components and therefore $\ta_\ho\cong\A$ is a theory.
The same equivalence induces an isomorphism of groups
$$
\Aut(f)\to\Aut(f_1)\x...\x\Aut(f_n),
$$
where $f:Y\to X$ is an 1-arrow in $\ta$ and $f_i=p_if:Y\to
X_n$. This fact together with the definition of a track extension
shows that $D$ is cartesian. Conversely, the above argument
actually shows that the functors
$$
\hog{Y,X}\to\hog{Y,X_1}\x...\x\hog{Y,X_n}
$$
induce bijections on connected components and induce isomorphisms
of corresponding automorphism groups. Hence these functors are
equivalences.
\end{proof}

Let $\A$ be a theory and let $D$ be  a natural system. As we saw
any object of the category $\Tracks(\A;D)$ is a track theory. We
now show that any morphism of $\Tracks(\A;D)$ carries lax products
to lax products. Indeed, let $F:\ta\to\ta'$ be a morphism of the
category $\Tracks(\A;D)$ and suppose all
$$
\hog{Y,X}_\ta\to\hog{Y,X_1}_\ta\x...\x\hog{Y,X_n}_\ta
$$
are equivalences of groupoids. We have to show that then the functors
$$
\hog{FY,FX}_{\ta'}\to\hog{FY,FX_1}_{\ta'}\x...\x\hog{FY,FX_n}_{\ta'},
$$
induced by $F$ are equivalences as well. But this is a consequence
of the following fact: For any $A,B\in\ta$ the functor
$F_{A,B}:\hog{A,B}_\ta\to\hog{FA,FB}_{\ta'}$ is an equivalence of
groupoids. To see the last assertion, it suffices to note that the
sets of components of both groupoids in question are canonically
isomorphic to $\Hom_\A(A,B)$ and $F$ is compatible with it and
also for any map $f:A\to B$ the groups $\Aut_\ta(f)$ and
$\Aut_{\ta'}(Ff)$ both are canonically isomorphic to $D_{qf}$.
Here $q:\ta_0\to\A$ is the canonical functor. According to Lemma
\ref{bololema} any morphism in $\Tracks(\A;D)$ is an equivalence
of theories.

\subsection{The main theorem}

We let
$\Str(\A;D)$ be the full subcategory of
$\Tracks(\A;D)$ whose objects are strong track theories.
There is an obvious functor
$$
\Str(\A;D)\to\Tracks(\A;D).
$$
Our Theorem \ref{theoria} shows that it yields a bijection on
the set of connected components.

The proof of Theorem \ref{theoria}
uses the following useful  construction.
Let
$$
0\to D\to \ta_1 \toto \ta_0 \xto p \bC \to 0
$$
be a linear track extension of a small category $\bC$ by a natural
system $D$. Here $D$ is any natural system on $\bC$. Suppose a
functor $f:\bE\to\ta_0$ is given which is identity on objects. We
assume that the composition $pf:\bE\to\bC$ is full. We construct
now a linear track extension
$$
0\to D\to\ta'_1\toto\ta'_0\xto{pf}\bC\to0
$$
with the property $\ta'_0=\bE$ and a morphism of linear track
extensions $\ta\to\ta'$ in $\Tracks(\bC;D)$. If $x,y:A\to B$ are
in $\bE$, then there exists a track $x\then y$  in $\ta'$ iff
$pf(x)=pf(y)$. If this holds, then we define the set
$\Hom_{\hog{A,B}_{\ta'}}(x,y)$ to be
$\Hom_{\hog{A,B}_\ta}(fx,fy)$. This defines a track extension
$\ta'$ which is denoted also by $f^!(\ta)$. The track functor
$f^!(\ta)\to\ta$ is identity on objects, on maps it is given by
$f$ and on tracks it is the inclusion.

\begin{The}\label{theoria}
Let $\A$ be a theory and let $D$ be a cartesian natural system on $\A$.
Then there exists a bijection
$$
\pi_0(\Str(\A;D))\cong H^3(\A;D).
$$
\end{The}

\begin{proof}
Let $\ta$ be an object of the category $\Str(\A;D)$. Then it can
be also considered as an object of $\Tracks(\A,\ta_0;D)$ and
therefore it defines an element in $H^3(\A,\ta_0;D)$ thanks to
Theorem \ref{fardobiti}. Then applying to this element the
boundary homomorphism $H^3(\A,\ta_0;D)\to H^3(\A;D)$ gives an
element in $H^3(\A;D)$. In this way we get a map
$$
\xi:\pi_0(\Str(\A;D))\to H^3(\A;D).
$$
We have to show that this map is a bijection. Take an $a\in
H^3(\A;D)$. There is a free theory $\F$ and a morphism of theories
$r:\F\to\A$ which is a full functor. Thanks to Theorem \ref{jib}
we have $H^i(\F;D)=0$ for all $i\ge2$. Therefore the connecting
homomorphism
$$
\partial:H^3(\A,\F;D)\to H^3(\A;D)
$$
is an isomorphism. Let $b=\partial\1(a)\in H^3(\A,\F;D)$ be the
element corresponding to $a$. Thanks to Theorem \ref{fardobiti}
the element $b$ defines a track extension $\ta$ with $\ta_0=\F$.
Thus $\ta$ is a strong track theory and hence $\xi$ is surjective.
It remains to show that $\xi$ is injective. Suppose
$\xi(\ta)=\xi(\ta')$. Since $p:\ta\to\A$ and $p':\ta'\to\A$ are
full morphisms of theories, one can lift the morphism of theories
$r:\F\to\A$ to morphisms $q:\F\to\ta_0$ and $q':\F\to\ta'_0$ of
theories, where $\F$ is a free theory and $r$ is surjective. Using
the $(\_ )^!$-construction one obtains the linear track extensions
$q^!(\ta)$ and ${q'}^!(\ta')$  together with morphisms of linear
track extensions $q^!(\ta)\to\ta$ and ${q'}^!(\ta')\to\ta'$. Now
both $q^!(\ta)$ and ${q'}^!(\ta')$ lie in $\Tracks(\A,\F;D)$, and
their classes in $ H^3(\A,\F;D)$ coincide with images of the
classes of $\ta$ and $\ta'$ under the homomorphisms $q^*:
H^3(\A,\ta_0;D)\to H^3(\A,\F;D)$ and ${q'}^*: H^3(\A,\ta'_0 ;D)\to
H^3(\A,\F;D)$ respectively. It follows from our assumptions that
these classes are the same and therefore $q^!(\ta)$ and
${q'}^!(\ta')$ are isomorphic in the groupoid $\Tracks(\A,\F;D)$.
Therefore we have the following diagram in $\Str(\A;D)$:
$$
\ta'\ot{q'}^!\ta'\ \cong \ q^!\ta \to \ta
$$
and hence the result.
\end{proof}

We are now in a position to prove our main result:

\begin{The}
Any abelian track theory is equivalent to a strong one. More
precisely if $\ta$ is an abelian track theory, then there exists a
strong abelian track theory $\ta'$, an abelian track theory
$\ta''$, and weak equivalences $\ta\leftarrow\ta''\rightarrow
\ta'$ as well as a lax equivalence $\ta'\to\ta$.
\end{The}

\begin{proof}
Let $\ta$ be an abelian track theory. Then the corresponding
homotopy category $\A:=\ta_\ho$ is a theory. Since any abelian
track category is part of a linear track extension, there is a
natural system $D$ on $\A$ such that $\ta\in\Tracks(\A;D)$. By
Proposition \ref{dekartulitrack} $D$ is a cartesian natural system
and therefore we can use Theorem \ref{theoria} to show that there
is an expected path in $\Tracks(\A;D)$ connecting $\ta$ to an
object of $\Str(\A;D)$. All maps in these diagrams are weak
equivalences thanks to Lemma \ref{bololema}. The last assertion
follows from Theorem \ref{laxfun}.
\end{proof}

\section*{Acknowledgments}

The second and third authors were partially supported by the grant
INTAS-99-00817 and by the RTN network HPRN-CT-2002-00287 ``K-theory and algebraic groups''.


\begin{thebibliography}{10}

\bibitem{sgaIV}
\emph{Th\'eorie des topos et cohomologie \'etale des sch\'emas. {T}ome 1:
  {T}h\'eorie des topos}, Springer-Verlag, Berlin, 1972, S\'eminaire de
  G\'eom\'etrie Alg\'ebrique du Bois-Marie 1963--1964 (SGA 4), Dirig\'e par M.
  Artin, A. Grothendieck, et J. L. Verdier. Avec la collaboration de N.
  Bourbaki, P. Deligne et B. Saint-Donat, Lecture Notes in Mathematics, Vol.
  269.

\bibitem{BaWe}
Michael Barr and Charles Wells, \emph{Toposes, triples and theories},
  Grundlehren der Mathematischen Wissenschaften [Fundamental Principles of
  Mathematical Sciences], vol. 278, Springer-Verlag, New York, 1985.

\bibitem{BE}
Hans-Joachim Baues, \emph{Combinatorial foundation of homology and homotopy},
  Springer Monographs in Mathematics, Springer-Verlag, Berlin, 1999,
  Applications to spaces, diagrams, transformation groups, compactifications,
  differential algebras, algebraic theories, simplicial objects, and
  resolutions.

\bibitem{Bsec}
Hans~Joachim Baues, \emph{The algebra of secondary cohomology operations},
  2001.

\bibitem{B4}
\bysame, \emph{The homotopy category of simply connected {$4$}-mainfolds},
  2003.

\bibitem{BD}
Hans~Joachim Baues and Winfried Dreckmann, \emph{The cohomology of homotopy
  categories and the general linear group}, $K$-Theory \textbf{3} (1989),
  no.~4, 307--338.

\bibitem{BJ1}
Hans-Joachim Baues and Mamuka Jibladze, \emph{Classification of abelian track
  categories}, $K$-Theory \textbf{25} (2002), no.~3, 299--311.

\bibitem{BTONKS}
Hans-Joachim Baues and Andy Tonks, \emph{On sum-normalised cohomology of
  categories, twisted homotopy pairs and universal toda brackets}, Quart. J.
  Math. \textbf{47} (1996), no.~188, 405--433.

\bibitem{BW}
Hans~Joachim Baues and G{\"u}nther Wirsching, \emph{Cohomology of small
  categories}, J. Pure Appl. Algebra \textbf{38} (1985), no.~2-3, 187--211.


\bibitem{BV}
J.~M. Boardman and R.~M. Vogt, \emph{Homotopy invariant algebraic structures on
  topological spaces}, Springer-Verlag, Berlin, 1973, Lecture Notes in
  Mathematics, Vol. 347. \MR{54 \#8623a}

\bibitem{fct}
John~W. Gray, \emph{Formal category theory: adjointness for {$2$}-categories},
  Springer-Verlag, Berlin, 1974, Lecture Notes in Mathematics, Vol. 391.

\bibitem{illusie}
Luc Illusie, \emph{Complexe cotangent et d\'eformations. {II}},
  Springer-Verlag, Berlin, 1972, Lecture Notes in Mathematics, Vol. 283.

\bibitem{JP2}
Mamuka Jibladze and Teimuraz Pirashvili, \emph{Cohomology of algebraic
  theories},  \textbf{137} (1991), no.~2, 253--296.

\bibitem{JPN}
\bysame, \emph{Linear extensions and nilpotence of maltsev theories}, 2000,
  SFB343 preprint 00-032, Universit\"at Bielefeld, p.~24.

\bibitem{working}
Saunders MacLane, \emph{Categories for the working mathematician},
  Springer-Verlag, New York, 1971, Graduate Texts in Mathematics, Vol. 5.

\bibitem{Mitchel}
Barry Mitchell, \emph{Rings with several objects}, Advances in Math. \textbf{8}
  (1972), 1--161.

\bibitem{P1}
T.~I. Pirashvili, \emph{Models for the homotopy theory and cohomology of small
  categories}, Soobshch. Akad. Nauk Gruzin. SSR \textbf{129} (1988), no.~2,
  261--264.

\bibitem{P2}
Teimuraz Pirashvili, \emph{Cohomology of small categories in homotopical
  algebra}, $K$-theory and homological algebra (Tbilisi, 1987--88), Lecture
  Notes in Math., vol. 1437, Springer, Berlin, 1990, pp.~268--302.

\bibitem{PW}
Teimuraz Pirashvili and Friedhelm Waldhausen, \emph{Mac {L}ane homology and
  topological {H}ochschild homology}, J. Pure Appl. Algebra \textbf{82} (1992),
  no.~1, 81--98.

\bibitem{R}
W.~H. Rowan, \emph{Enveloping ringoids}, Algebra Universalis \textbf{35}
  (1996), no.~2, 202--229.

\bibitem{Schwede}
Stefan Schwede, \emph{Stable homotopy of algebraic theories}, Topology
  \textbf{40} (2001), no.~1, 1--41.

\end{thebibliography}
\end{document}